\documentclass[12pt]{amsart}

\usepackage{graphicx, comment}
\usepackage{hyperref}

\usepackage{biblatex}
\usepackage[normalem]{ulem}
\usepackage{tikz}
\usepackage[normalem]{ulem}
\usepackage{pgfplots}
\usepackage{multicol}
\usepackage{standalone}[mode=buildnew]
\usetikzlibrary{calc,3d}
\usetikzlibrary{arrows.meta,calc,decorations.markings,math,arrows.meta, intersections}
\usepgfplotslibrary{fillbetween}
\usepackage{xcolor}

\usepackage{enumitem}

\usepackage{algorithm2e}

\usepackage[margin=3cm]{geometry}
\usepackage{caption}
\usepackage{subcaption}
\SetKw{Or}{\hspace{\algoskipindent}\itshape or\;}
\SetKw{And}{\hspace{\algoskipindent}\itshape and\;}
\SetKwBlock{Condition}{}{}

\theoremstyle{definition}

\theoremstyle{definition}

\theoremstyle{remark}

\numberwithin{equation}{section}

\newcounter{problem}[section]
\newcommand\be{\begin{equation*}}
\newcommand\ee{\end{equation*}}
\newcommand\bea{\begin{eqnarray*}}
\newcommand\eea{\end{eqnarray*}}
\newcommand\bi{\begin{itemize}}
\newcommand\ei{\end{itemize}}
\newcommand\bean{\begin{equation}}
\newcommand\eean{\end{equation}}
\newcommand\ben{\begin{enumerate}}
\newcommand\een{\end{enumerate}}
\newtheorem{thm}{Theorem}[section]

\newtheorem{cor}[thm]{Corollary}
\newtheorem{lem}[thm]{Lemma}
\newtheorem{prop}[thm]{Proposition}

\newtheorem{defi}[thm]{Definition}

{\refstepcounter{problem}}

\let\oldnl\nl% Store \nl in \oldnl
\newcommand{\nonl}{\renewcommand{\nl}{\let\nl\oldnl}}% Remove line number for one line

\SetKw{KwBy}{by}

\newtheorem*{lem*}{Lemma}

\newtheorem{innercustomthm}{Theorem}
\newenvironment{customthm}[1]
  {\renewcommand\theinnercustomthm{#1}\innercustomthm}
  {\endinnercustomthm}

\newcommand\floor[1]{\left\lfloor#1\right\rfloor}

\newcommand{\R}{\ensuremath{\mathbb{R}}}
\newcommand{\C}{\ensuremath{\mathbb{C}}}
\newcommand{\Z}{\ensuremath{\mathbb{Z}}}
\newcommand{\Q}{\mathbb{Q}}
\newcommand{\N}{\mathbb{N}}

\newcommand{\norm}[1]{\left|\left|#1\right|\right|}
\newcommand{\m}{\textbf{m}}
\newcommand{\w}{\textbf{w}}

\makeatletter
\newcommand\subsubsubsection{\@startsection{paragraph}{4}{\z@}{-2.5ex\@plus -1ex \@minus -.25ex}{1.25ex \@plus .25ex}{\normalfont\normalsize\bfseries}}
\newcommand\subsubsubsubsection{\@startsection{subparagraph}{5}{\z@}{-2.5ex\@plus -1ex \@minus -.25ex}{1.25ex \@plus .25ex}{\normalfont\normalsize\bfseries}}
\makeatother

\title[RF-Values for Embeddings of 4-D Ellipsoids into Almost-Cubes]{Rigid-Flexible Values for Embeddings of Four-Dimensional Ellipsoids into Almost-Cubes}
\author{Cory H. Colbert}
\address{Washington and Lee University, Lexington, VA 24450}
\email{ccolbert@wlu.edu}

\author{Andrew Lee}
\address{St. Thomas Aquinas College, Sparkill, NY 10976}
\email{alee@stac.edu}


\begin{document}
\begin{abstract}
%In previous work of Cristofaro-Gardiner, Frenkel, and Schlenk,  was determined for all integers $b \geq 2.$
%As in McDuff's work, Cristofaro-Gardiner, Frenkel, and Schlenk found staircases associated with the embedding function.  More recently, Usher's work shows that the appearance of infinite staircases as $b$ varies is hard to predict. The intricate structure found there suggests that determining the entirety of the graph of $c_b(a)$ for all $b$ is intractable.

We consider the embedding function $c_b(a)$ describing the problem of symplectically embedding an ellipsoid $E(1,a)$ into the smallest possible scaling by $\lambda>1$ of the polydisc $P(1,b)$. In particular, we calculate rigid-flexible values, i.e. the minimum $a$ such that for $E(1,a')$ with $a'>a$, the embedding problem is determined only by volume. For $1<b<2$ we find that these values vary piecewise smoothly outside the discrete set $b\in\left(\frac{n+1}{n}\right)^2$. As Jin and Lee analyze packing stability in the caes $b>2$, our results complete the story outside of a discrete set.
\end{abstract}

\maketitle
\thispagestyle{empty}
%\newpage

\tableofcontents

\section{Introduction, Statement of Results}

The problem of embedding one symplectic manifold into another touches on a wide variety of topics in symplectic geometry, and in this work we focus in particular on embeddings of ellipsoids $E(a,b)$ into polydiscs $P(a,b)$ . Here, a polydisc
\be
P(a,b) := \left\{(z_1,z_2) \in \C^2 \,| \, \pi \|z_1\|^2 < a, \pi \|z_2\|^2 < b \right\} \nonumber
\ee
is the 4-dimensional open symplectic manifold $B_2(a)\times B_2(b)\subset \C^2$, where each factor is a 2-disc of fixed radius centered at $0\in\C$. Similarly the ellipsoid $E(a,b)$ is given by
\be
E(a,b) := \left\{(z_1, z_2)\in\C^2\,|\, \frac{\pi\norm{z_1}^2}{a}+\frac{\pi\norm{z_2}^2}{b} < 1\right\}.\nonumber
\ee
% McDuff and Schlenk \cite{McSc} consider the problem of embedding a symplectic ellipsoid $E(1,a)$ into a ball $B(r)$: for fixed $a$, what is the smallest $r$ such that $E(1,a)\hookrightarrow B(r)$ symplectically? One has the immediate necessary condition that the volume of $E(1,a)$ must be no greater than that of $B(r)$ since the embedding is symplectic, but in fact this is far from sufficient. In \cite{McSc}, the authors determine the function $c(a)$ whose output is the minimal $r$ guaranteeing such an embedding, and show that the problems fluctuate between the flexibility seen in volume-preserving geometry and the rigidity often seen in symplectic geometry.
% In particular, \cite{McSc} find that if $a\geq \left(\frac{17}{6}\right)^2$, then $c(a) = \sqrt{a}$. %This bound is sharp, and records the $a$-value past which symplectic obstructions vanish.
We study the problem of finding embeddings of an ellipsoid into a polydisc. This information is encoded by a function $c_b(a)$ whose value at $a$ is the smallest $\lambda$ such that
\begin{equation}
E(1,a)\hookrightarrow P(\lambda, \lambda b).\label{embeddingproblem}
\end{equation}
More precisely,
\be
c_b(a):= \inf \{ \lambda > 0 \,|\, E(1,a)\hookrightarrow P(\lambda, \lambda b)\}.\nonumber
\ee
Earlier work on embeddings of ellipsoids into polydiscs, e.g. Frenkel-M\"uller and Cristofaro-Gardiner-Frenkel-Schlenk \cite{CFS, FM}, considered a fixed value of $b$ and calculated $c_b(a)$ in its entirety for all $a$. For example, $c_1(a)$ contains a so-called ``infinite staircase"a convergent sequence of $a$-values $a_n\to a_{\infty}$ such that $c_b(a)$ is non-smooth at the $a_i$ and alternates between constant and linear increasing on each $[a_i, a_{i+1}]$. %\textcolor{blue}{We call such an interval a \textit{step}, and the $a\in [a_i, a_{i+1}]$ where the function becomes constant is the \textit{center} of the step.}

In addition, they find that for sufficiently large values of $a$, $c_b(a)$ coincides with the volume constraint $\sqrt{\frac{a}{2b}}$: for sufficiently elongated ellipsoids the only obstruction to the existence of symplectic embeddings is the volume constraint. In the work of both Cristofaro-Gardiner-Frenkel-Schlenk and Usher \cite{CFS, Ush} both the source and target in the embedding problem were allowed to become elongated, with the goal of analyzing the resulting changes in the embedding function $c_b(a)$.

Varying $b$ uncovers delicate structure in the the $b$-direction for functions $c_b(a)$. In Cristofaro-Gardiner-Frenkel-Schlenk \cite{CFS} the authors find that for $b\in\N$ with $b\geq 2$, the function $c_b(a)$ contains no infinite staircases and can be computed explicitly for all $a$. On the other hand, Usher's work in \cite{Ush} provides examples of irrational $b$ (some arbitrarily large and others arbitrarily close to integers) such that for each such $b$, $c_b(a)$ contains an infinite staircase. In the same work, Usher also computes the embedding functions $c_b(a)$ for $b\in [1,\infty)$ and $a\in [1,3+\sqrt{2}]$.

It is also of interest, then, if and when this complexity disappears: i.e. whether at some point, $c_b(a)$ is guaranteed to coincide with the volume constraint. At such a point, the embedding problem becomes \emph{flexible}, in the sense that the answer only depends on volume, and obstructions coming from \emph{rigid} geometric objects (in the form of holomorphic curves) no longer appear. For this reason, we focus on a particular portion of the graph of $c_b(a)$ and define the \emph{rigid-flexible value}, or $RF$-value, as
\begin{equation}
RF(b): = \inf\left\{A\left| c_b(a) = \sqrt{\frac{a}{2b}}\,\,\mbox{for all}\,\, a\geq A\right.\right\}.\label{rfdefn}
\end{equation}
%Background on what's known about RF(b) and its properties
The fact that $RF$ is even a well-defined real-valued function is nontrivial. Indeed, in \cite[Theorem 1.1]{BPT}, the authors show that in (\ref{embeddingproblem}), we must have $c_b(a) = \sqrt{\frac{a}{2b}}$ for $a \geq \frac{9(b+1)^2}{2b}$, so the set in (\ref{rfdefn}) contains a positive real number.

In the search for rigidity in these embedding problems, the set of so-called \textit{exceptional classes} plays a special role. These are classes in the second homology of $\mathbb{C}P^2$ blown up in multiple points: specifically, those represented by embedded holomorphic spheres of self-intersection -1. The symplectic blow-up construction described in \ Work of Li-Liu in \cite{LiLiu} and Li-Li in \cite{LiLi} led to a combinatorial characterization of these classes, unlocking computationally tractable methods for calculating $c_b(a)$ in a range of different embedding problems. Each such class represents a potential obstruction to fully filling the polydisc with an ellipsoid, and the $RF$-value is the precise ellipsoid length $a$ such that no more such obstructions appear.

The variation of $c_b(a)$ with respect to $b$ uncovered by these methods makes computing $RF$ a nontrivial task. For example, Cristofaro-Gardiner-Frenkel-Schlenk \cite[Section 3.6]{CFS} shows that $RF(2) = 8\frac{1}{36}$ given by the exceptional class $(6,3;3,2^{\times 7})$, shown as the only filled dot in Figure \ref{rffigure:bonetotwo}. Frenkel and M\"uller in \cite[Section 7.3]{FM} show that $RF(1) = 7\frac{1}{36}$ as determined by the exceptional class $(4,4; 3, 2^{\times 6})$. The main result in Jin and Lee's work\cite{JL} is entirely devoted to the behavior of $RF$, where the authors show that $RF(b)$ for $b > 2$ is piecewise smooth and right continuous on the interval $[2, +\infty)$, with an infinite discrete set of discontinuities at $\{n+2-\sqrt{2n+3}: n\in\mathbb N\}$. In that work, the authors also show that $RF$ is not continuous at $b = 1$ by demonstrating that $\lim_{m \to \infty} RF\left(1 + \frac{1}{m}\right) = 8,$ contrasting with the result of Frenkel and M\"uller \cite{FM} at $b=1$ cited above. These results narrow our study of $RF$ to the interval $(1, 2),$ which is where our work in this paper is focused. We term the polydiscs of the form $P(\lambda, \lambda b)$ for $b\in (1,2)$ ``almost-cubes" following Frenkel and M\"uller, where cubes are polydiscs of the form $P(\lambda,\lambda)$. We determine $RF$ on $(1,2)$ precisely outside the discrete set of ratios of consecutive perfect squares $\left\{\frac{(n+1)^2}{n^2}: n \in \mathbb N\right\}.$ In this work, we consider the classes
\begin{eqnarray}
S_k &:=& (2k^2 + k, 2k^2 - k; (k^2)^{\times 7}, k^2-1),\label{skclasses}\\
T_k &:=& ((2k + 1)(k + 1), k(2k + 1); (k^2 + k + 1), (k^2 + k)^{\times 7}).\label{tkclasses}
\end{eqnarray}
where the $T_k$ appeared in Jin and Lee's work\cite{JL}.

% Previous joint work \cite{JL} by the second author computed $RF(b)$ for $b>2$, finding that it was piecewise smooth in $b$ away from an infinite discrete set of discontinuities. \textcolor{red}{These discontinuities occur at values of the function $u(n) = n + 2 - \sqrt{2n + 3}$ for $n\in \mathbb{N}$. Moreover, for $b\geq 2$, $RF(b)$ is right continuous.}

% \cite{JL} computed $\lim_{m\to\infty}RF\left(\frac{m+1}{m}\right) = 8$, which is consistent with the fact from \cite{CFS} that $RF(2) = 8\frac{1}{36}$, determined by the exceptional class $(6,3;3,2^{\times 7})$. However, \cite{FM} show that $RF(1)=7\frac{1}{32}$, determined by the exceptional class $(4,4;3,2^{\times 6})$. This fact, together with the conclusion of \cite{JL} on the discrete set $\left{\frac{n+1}{n}\right}$, further established that $RF(b)$ was not continuous at $b=1$.

% In this paper we determine $RF$ for $1<b<2$ outside a discrete set. We term these polydiscs $P(\lambda, \lambda\cdot b)$ ``almost-cubes" following \cite{FM}, where cubes are polydiscs $P(\lambda,\lambda)$. In this region we again see an infinite set of discontinuities, this time accumulating at 1, but away from these points $RF(b)$ is smooth.

For each $n \ge 1,$ define the open intervals

\be
I_n := \left(\left(\frac{n+2}{n+1}\right)^2,\left(\frac{n+1}{n}\right)^2\right),
\ee

and let $\mu_b(C)(a)$ be the obstruction function associated to the class $C$ as in Theorem \ref{classcondition}, using the notation of \cite{CFS, FM, McSc}.

\begin{thm}\label{mainThm}
\textit{
Let $b \in \left(1, 2 \right)\setminus\left\{\frac{(n+1)^2}{n^2}: n \in \mathbb N\right\}.$ Then
\begin{equation}
RF(b) =
\begin{cases}
    2b(\mu_b(S_{k+1})(8))^2 & \mbox{if $b \in I_{2k+1}$} \\
    2b(\mu_b(T_k)(8))^2 & \mbox{if $b \in I_{2k}$} \\
   \end{cases}.
\end{equation}
}
\end{thm}
\begin{figure}
\centering
\scalebox{0.9}{
\begin{subfigure}[t]{.5\textwidth}
  \captionsetup{width=.9\linewidth}
  %\begin{questions}
\begin{tikzpicture}
\begin{axis}[
            xtick={1.5,2},
            xtick style = {at={(axis description cs:0,0.7)}, anchor=north},
            ytick={8},
            xmin=1.2,
            xmax=2,
            ymin=7.98,
            ymax=8.03]
%At 25/16-0.0001=1.5624 class at 8 is 15,10;7,6,6,6,6,6,6,6 = T_2
%At 25/16+0.0001=1.5626 class at 8 is 10,6: 4,4,4,4,4,4,4,3 = S_2
%At 36/25-0.0001=1.4399 class at 8 is 21,15;9,9,9,9,9,9,9,8 = T_3
%
\addplot[gray, thin, domain=1.2:2]{8};
%Add in plot pieces
\addplot[gray!80, thick, domain=1.7777777777777777:2] { 2*x*((8*1^2+8*1+1)/(2*1^2+3*1+1+2*x*1^2+x*1))^2 } node[black, above,pos=0.7, font=\tiny] {{$T_1$}};
%\addplot[domain=1.7777777777777777:2, gray!50] { 2*x*((8*1^2-1)/(2*1^2+1+2*x*1^2-x*1))^2 } node[below,pos=0.5, font=\tiny] {{$T_1$}};
\addplot[gray!80, thick, domain=1.5625:1.7777777777777777] { 2*x*((8*2^2-1)/(2*2^2+2+2*x*2^2-x*2))^2 } node[black, above,pos=0.5, font=\tiny] {{$S_2$}};
%\addplot[domain=1.5625:1.7777777777777777, gray!50] { 2*x*((8*2^2+8*2+1)/(2*2^2+3*2+1+2*x*2^2+x*2))^2 } node[below left,pos=0.9, font=\tiny] {{$T_2$}};
\addplot[gray!80, thick, domain=1.44:1.5625] { 2*x*((8*2^2+8*2+1)/(2*2^2+3*2+1+2*x*2^2+x*2))^2 } node[black, above,pos=0.5, font=\tiny] {{$T_2$}};
%\addplot[domain=1.44:1.5625, gray!50] { 2*x*((8*2^2-1)/(2*2^2+2+2*x*2^2-x*2))^2 } node[below right,pos=0.2, font=\tiny] {{$S_2$}};
\addplot[gray!80, thick, domain=1.3611111111111112:1.44] { 2*x*((8*3^2-1)/(2*3^2+3+2*x*3^2-x*3))^2 } node[black, above,pos=0.5, font=\tiny] {{$S_3$}};
\addplot[gray!40, thick, domain=1.3611111111111112:1.44, gray!50] { 2*x*((8*3^2+8*3+1)/(2*3^2+3*3+1+2*x*3^2+x*3))^2 } node[gray!40, below,pos=0.5, font=\tiny] {{$T_3$}};%below RF
\addplot[gray!80, thick, domain=1.3061224489795917:1.3611111111111112] { 2*x*((8*3^2+8*3+1)/(2*3^2+3*3+1+2*x*3^2+x*3))^2 } node[black, above,pos=0.5, font=\tiny] {{$T_3$}};
\addplot[domain=1.3061224489795917:1.3611111111111112, gray!50] { 2*x*((8*3^2-1)/(2*3^2+3+2*x*3^2-x*3))^2 } node[below,pos=0.5, font=\tiny] {{$S_3$}};
\addplot[gray!80, thick, domain=1.265625:1.3061224489795917, samples=5] { 2*x*((8*4^2-1)/(2*4^2+4+2*x*4^2-x*4))^2 } node[black, above,pos=0.5, font=\tiny] {{$S_4$}};
\addplot[gray!50, thick, domain=1.265625:1.3061224489795917, gray!50, samples=4] { 2*x*((8*4^2+8*4+1)/(2*4^2+3*4+1+2*x*4^2+x*4))^2 } node[below,pos=0.5, font=\tiny] {{$T_4$}};
\addplot[domain=1.2345679012345678:1.265625, samples=3] { 2*x*((8*4^2+8*4+1)/(2*4^2+3*4+1+2*x*4^2+x*4))^2 } node[above,pos=0.5, font=\tiny] {{$T_4$}};
\addplot[domain=1.2345679012345678:1.265625, gray!50, samples=3] { 2*x*((8*4^2-1)/(2*4^2+4+2*x*4^2-x*4))^2 } node[below,pos=0.5, font=\tiny] {{$S_4$}};

%Add in endpoints
\addplot[mark=*, fill=black, only marks, mark size=1pt] coordinates {(2, 8.027777777)};
\addplot[mark=*, fill=white, only marks, mark size=1pt] coordinates {(1.7777777777777777,7.945854639471415)};
\addplot[mark=*, fill=white, only marks, mark size=1pt] coordinates {(1.5625,8.000000000000002)};
\addplot[mark=*, fill=white, only marks, mark size=1pt] coordinates {(1.44,8.0)};
\addplot[mark=*, fill=white, only marks, mark size=1pt] coordinates {(1.44,7.989014762709816)};
\addplot[mark=*, fill=white, only marks, mark size=1pt] coordinates {(1.3611111111111112,7.999999999999999)};
\addplot[mark=*, fill=white, only marks, mark size=1pt] coordinates {(1.3061224489795917,8.000000000000002)};
\addplot[mark=*, fill=white, only marks, mark size=1pt] coordinates {(1.3061224489795917,7.996451934847931)};
\addplot[mark=*, fill=white, only marks, mark size=1pt] coordinates {(1.265625,8.0)};
\addplot[mark=*, fill=white, only marks, mark size=1pt] coordinates {(1.2345679012345678,8.0)};
%Tks
\addplot[mark=*, fill=white, only marks, mark size=1pt] coordinates {(1.44,7.965702075926983)};
\addplot[mark=*, fill=white, only marks, mark size=1pt] coordinates {(1.5625,8.000000000000002)};
\addplot[mark=*, fill=white, only marks, mark size=1pt] coordinates {(1.7777777777777777,8.000000000000004)};
\addplot[mark=*, fill=white, only marks, mark size=1pt] coordinates {(1.3061224489795917,7.991957778844478)};
\addplot[mark=*, fill=white, only marks, mark size=1pt] coordinates {(1.3611111111111112,8.000000000000002)};
\addplot[mark=*, fill=white, only marks, mark size=1pt] coordinates {(1.44,7.999999999999998)};
\addplot[mark=*, fill=white, only marks, mark size=1pt] coordinates {(1.2345679012345678,7.997201084601139)};
\addplot[mark=*, fill=white, only marks, mark size=1pt] coordinates {(1.265625,8.0)};
\addplot[mark=*, fill=white, only marks, mark size=1pt] coordinates {(1.3061224489795917,8.000000000000004)};
%\legend{y=x+3}
\end{axis}
\end{tikzpicture}
  \caption{$RF(b)$ for $b\in[1.2, 2]$. Labels indicate whether $\mu(S_k)$ or $\mu(T_k)$ determines $RF$ on the given open interval. The graphs of some obstruction functions are continued in gray to exhibit the exchange between $S_k$ and $T_k$. The values of $RF$ at the open dots are unknown.}
  \label{rffigure:bonetotwo}
\end{subfigure}
}%\hfill
\scalebox{0.9}{
\begin{subfigure}[t]{.5\textwidth}
\captionsetup{width=.9\linewidth}
  \begin{tikzpicture}
  \begin{axis}[
%    xlabel={\tiny{$b$, Polydisc length}}, % default put x on x-axis
%    ylabel={\tiny{$a$, Ellipsoid length}},  % default put y on y-axis
%    x label style = {at={(axis description cs:0.45, -0.01)},anchor=north},
%    y label style = {at={(axis description cs:0,0.4)}, anchor=south},
    samples=100,
    ytick={8,12},
    yticklabels={$8$,$12$},
    xtick={1,2,3,4},
    xticklabels={$1$,$2$,$3$,$4$},
    domain=1.25:6,
    xmin=1.26,xmax=4,
    ymin=7.5,ymax=13.5,
    grid=both, grid style={line width=.1pt}
  ]
%  \node [anchor=west] at (axis cs:  1.25, 5.5) {\tiny{$3+2\sqrt{2}$}};
  %\node at (axis cs:  2.5,  6.2) {\small{\textcolor{magenta}{Infinite staircases live here!}}};
%  \node [anchor=south] at (axis cs: 2, 1) {2};
  %\node at (1.75,8.5) {\tiny{\textcolor{red}{RF}}};
  %\node at (axis cs:  2.45, 8.46) {\tiny{No staircases here}};
  %\node at (axis cs:  2.75, 3.5) {\small{Fully understood; no staircases past $b=1$}};
%  \node at (axis cs:  4, .5) {\tiny{}};
%  \node at (axis cs:  4, 2.5) {\tiny{}};
%  \node at (axis cs:  3.5, 3.5) {\tiny{Integer $\beta$ only needs $E_n, F_n$}};
%  \node at (axis cs:  2.7, 3) {\tiny{when $\beta\geq 2$}};
%  \node at (axis cs:  3.4,2.5) {\tiny{$\exists$ infinite staircases for inf.}};
%  \node at (axis cs:  2.9,2) {\tiny{ many irrational $\beta$}};
%  \node at (axis cs:  1.5,2) {\tiny{}};
%  \node [anchor=south west] at (axis cs: 1.25,1) {$(1,1)$};
  % \addplot[name path=NSlesstwo, green, thick, mark=none, domain=1.265:5] {((2*x+2)/sqrt(2*x)+sqrt((2*x+2)*(2*x+2)/(2*x)-4))*((2*x+2)/sqrt(2*x)+sqrt((2*x+2)*(2*x+2)/(2*x)-4))/4}; %No staircases line
%   lincoeff=(2*b+2)/np.sqrt(2*b)
%   root=(lincoeff+np.sqrt(lincoeff*lincoeff-4))/(2)
%   limit=root*root
%  \addplot[name path=betatwo, gray, dotted, mark=none, domain=0:5] {2.25};
%  \addplot[name path=NS, green, thick, mark=none, domain=3.63:4.9] {x*x-3*x};
%  \addplot[name path=threeplus, blue, thick, mark=none] coordinates {(1.25, 5.8) (5, 5.8)};
%    \addplot[name path=notdef, black, thick, mark=none] coordinates {(1.25, 1) (5, 1)};
  \addplot[name path=RFstep1, black, thick, mark=none, domain=2:2.6834] {2*x*(9/(x+4))*(9/(x+4))};
  \addplot[name path=RFstep2, black, thick, mark=none, domain=2.6834:3.3945] {2*x*(11/(x+5))*(11/(x+5))};
  \addplot[name path=RFstep3, black, thick, mark=none, domain=3.3945:4.127] {2*x*(13/(x+6))*(13/(x+6))};
  \addplot[name path=RFstep4, black, thick, mark=none, domain=4.127:4.877] {2*x*(15/(x+7))*(15/(x+7))};
%  \addplot[name path=RFunknown, black, dotted, mark=none, domain=4:4.25] {x*x-4*x};
% \addplot fill between[
%     of = NS and RF,
%     soft clip={domain=3:5},
%     every even segment/.style  = {gray,opacity=.4}
%   ];
% \addplot fill between[
%     of = NSlesstwo and RF,
%     soft clip={domain=2:4.5},
%     every even segment/.style  = {gray,opacity=.4}
%   ];
%  \draw[thick,dashed,brown] (axis cs:0,3) -- (axis cs:0,2);

%Add in fill for each plot piece
\addplot[name path= RF1, gray, thick, mark=none, domain=1.7777777777777777:2] {2*x*((8*1^2+8*1+1)/(2*1^2+3*1+1+2*x*1^2+x*1))^2} ;
% \addplot[red, thick, mark=none, domain=1.7777777777777777:2, gray!50] { 2*x*((8*1^2-1)/(2*1^2+1+2*x*1^2-x*1))^2 };
% \addplot fill between[
%     of = NSlesstwo and RF1,
%     soft clip={domain=1.777777777777:2},
%     every even segment/.style  = {gray,opacity=.4}
%   ];
\addplot[name path = RF2, gray, thick, mark=none, domain=1.5625:1.7777777777777777] { 2*x*((8*2^2-1)/(2*2^2+2+2*x*2^2-x*2))^2 };
% \addplot fill between[
%     of = NSlesstwo and RF2,
%     soft clip={domain=1.5625:1.7777777777777777},
%     every even segment/.style  = {gray,opacity=.4}
%   ];
%\addplot[red, thick, mark=none, domain=1.5625:1.7777777777777777, gray!50] { 2*x*((8*2^2+8*2+1)/(2*2^2+3*2+1+2*x*2^2+x*2))^2 };
\addplot[name path = RF3, gray, thick, mark=none, domain=1.44:1.5625] { 2*x*((8*2^2+8*2+1)/(2*2^2+3*2+1+2*x*2^2+x*2))^2 };
% \addplot fill between[
%     of = NSlesstwo and RF3,
%     soft clip={domain=1.44:1.5625},
%     every even segment/.style  = {gray,opacity=.4}
%   ];
%\addplot[red, thick, mark=none, domain=1.44:1.5625, gray!50] { 2*x*((8*2^2-1)/(2*2^2+2+2*x*2^2-x*2))^2 };
\addplot[name path = RF4, gray, thick, mark=none, domain=1.3611111111111112:1.44] { 2*x*((8*3^2-1)/(2*3^2+3+2*x*3^2-x*3))^2 };
% \addplot fill between[
%     of = NSlesstwo and RF4,
%     soft clip={domain=1.3611111111111112:1.44},
%     every even segment/.style  = {gray,opacity=.4}
%   ];
%\addplot[red, thick, mark=none, domain=1.3611111111111112:1.44, gray!50] { 2*x*((8*3^2+8*3+1)/(2*3^2+3*3+1+2*x*3^2+x*3))^2 };
\addplot[name path = RF5, gray, thick, mark=none, domain=1.3061224489795917:1.3611111111111112] { 2*x*((8*3^2+8*3+1)/(2*3^2+3*3+1+2*x*3^2+x*3))^2 };
% \addplot fill between[
%     of = NSlesstwo and RF5,
%     soft clip={domain=1.3061224489795917:1.3611111111111112},
%     every even segment/.style  = {gray,opacity=.4}
%   ];
%\addplot[red, thick, mark=none, domain=1.3061224489795917:1.3611111111111112, gray!50] { 2*x*((8*3^2-1)/(2*3^2+3+2*x*3^2-x*3))^2 };
\addplot[name path = RF6, gray, thick, mark=none, domain=1.265625:1.3061224489795917, samples=5] { 2*x*((8*4^2-1)/(2*4^2+4+2*x*4^2-x*4))^2 };
% \addplot fill between[
%     of = NSlesstwo and RF6,
%     soft clip={domain=1.265625:1.3061224489795917},
%     every even segment/.style  = {gray,opacity=.4}
%   ];
%\addplot[red, thick, mark=none, domain=1.265625:1.3061224489795917, gray!50, samples=4] { 2*x*((8*4^2+8*4+1)/(2*4^2+3*4+1+2*x*4^2+x*4))^2 };
%\addplot[red, thick, mark=none, domain=1:1.265625, samples=3] { 2*x*((8*4^2+8*4+1)/(2*4^2+3*4+1+2*x*4^2+x*4))^2 };
%\addplot[red, thick, mark=none, domain=1.2345679012345678:1.265625, gray!50, samples=3] { 2*x*((8*4^2-1)/(2*4^2+4+2*x*4^2-x*4))^2 };

% %Add in endpoints
\addplot[mark=*, fill=black, only marks, mark size=1pt] coordinates {(2, 8.027777777)};
%\addplot[mark=*, fill=white, only marks, mark size=1pt] coordinates {(1.7777777777777777,7.945854639471415)};
%\addplot[mark=*, fill=white, only marks, mark size=1pt] coordinates {(1.5625,8.000000000000002)};
%\addplot[mark=*, fill=white, only marks, mark size=1pt] coordinates {(1.44,8.0)};
%\addplot[mark=*, fill=white, only marks, mark size=1pt] coordinates {(1.44,7.989014762709816)};
%\addplot[mark=*, fill=white, only marks, mark size=1pt] coordinates {(1.3611111111111112,7.999999999999999)};
%\addplot[mark=*, fill=white, only marks, mark size=1pt] coordinates {(1.3061224489795917,8.000000000000002)};
%\addplot[mark=*, fill=white, only marks, mark size=1pt] coordinates {(1.3061224489795917,7.996451934847931)};
%\addplot[mark=*, fill=white, only marks, mark size=1pt] coordinates {(1.265625,8.0)};
%\addplot[mark=*, fill=white, only marks, mark size=1pt] coordinates {(1.2345679012345678,8.0)};
%Tks
%\addplot[mark=*, fill=white, only marks, mark size=1pt] coordinates {(1.44,7.965702075926983)};
\addplot[mark=*, fill=white, only marks, mark size=1pt] coordinates {(1.5625,8.000000000000002)};
\addplot[mark=*, fill=white, only marks, mark size=1pt] coordinates {(1.7777777777777777,8.000000000000004)};
\addplot[mark=*, fill=white, only marks, mark size=1pt] coordinates {(1.3061224489795917,7.991957778844478)};
\addplot[mark=*, fill=white, only marks, mark size=1pt] coordinates {(1.3611111111111112,8.000000000000002)};
\addplot[mark=*, fill=white, only marks, mark size=1pt] coordinates {(1.44,7.999999999999998)};
\addplot[mark=*, fill=white, only marks, mark size=1pt] coordinates {(1.26,7.997201084601139)};
\addplot[mark=*, fill=white, only marks, mark size=1pt] coordinates {(2,9.0)};
\addplot[mark=*, fill=white, only marks, mark size=1pt] coordinates {(1.3061224489795917,8.000000000000004)};

%points for RFsteps
\addplot[mark=*, fill=black, only marks,    mark size=1pt] coordinates {(2.68338,11)};
\addplot[mark=*, fill=black, only marks, mark size=1pt] coordinates {(3.39445,13)};
\addplot[mark=*, fill=black, only marks, mark size=1pt] coordinates {(4.12702,15)};

\addplot[mark=*, fill=white, only marks, mark size=1pt] coordinates {(2.683,9.73206)};
\addplot[mark=*, fill=white, only marks, mark size=1pt] coordinates {(3.39445,11.65736)};
\addplot[mark=*, fill=white, only marks, mark size=1pt] coordinates {(4.12702,13.6016)};

\end{axis}
\end{tikzpicture}
  \caption{$RF(b)$ for $b\in [1,4]$. The leftmost portion of the figure (in gray) is depicted in Figure \ref{rffigure:bonetotwo}. The values of $RF(b)$ when $b\in \mathbb{N}$ are established in \cite{FM} for $b=1$ and $b\geq 2$ in \cite{CFS}. The values for $b\in\R$ where in addition $RF(b)$ is shown to be right continuous.}
  \label{rffigure:bonetothirteen}
\end{subfigure}
}
\caption{Graphs of $RF(b)$, with the aspect ratio $b$ of the polydisc $P(\lambda, \lambda b)$ on the horizontal axis.}
\end{figure}
As an example, for $b\in\left(\frac{25}{16}, \frac{16}{9}\right) = I_3,$ we have $3=2k+1$ so we consider the class $S_2 = (10, 6; 4^{\times 7}, 3)$, which will obstruct the existence of a symplectic embedding of $E(1,8)$ into a scaling of $P(1,b)$ of the same volume.. By Theorem \ref{mainThm} $RF(b)$ for $b\in I_3$ is given by $2b\left(\frac{31}{10+6b}\right)^2$. To see how the relevant class changes with $b$, note that for $b\in\left(\frac{36}{25},\frac{25}{16}\right)=I_4$, we have $2k=4$ so Theorem \ref{mainThm} points us to the $T_k$-classes with $k=2$. In this interval the class $T_2 = (15,10;7, 6^{\times 7})$ provides an obstruction to embedding $E(1,8)$ symplectically into a scaling of $P(1,b)$ with the same volume.

Moreover, by Theorem \ref{mainThm} $RF(b)$ for $b\in I_4$ is given by $2b\left(\frac{49}{15+10b}\right)^2$. Continuing in this way we can plot $RF(b)$ for each $b\in \left(\frac{(n+2)^2}{(n+1)^2},\frac{(n+1)^2}{n^2}\right)$. This is the source of Figures \ref{rffigure:bonetothirteen} and \ref{rffigure:bonetotwo}.

Notably, our $T_1$ coincides with $F_2$ of Cristofaro-Gardiner-Frenkel-Schlenk, which provides a nice transition from their work to ours. The family of classes $F_b$ defined therein contains classes with arbitrarily many entries and produces obstructions for $b\in \mathbb{N}^{\geq 2}$, while our $S_k$ and $T_k$ all have exactly 10 entries for any $k$ and produce obstructions for $1<b<2$.

% \begin{figure}\label{rfgeography}
% \include{embeddinggeography}
% \caption{Geography of the embedding problem $E(1,a)\hookrightarrow P(\lambda,\lambda b)$ for $1\leq a \leq 9$ and $1\leq b \leq ??$.}
% \end{figure}

On the other hand, while $RF(b)$ for $b\geq 2$ is right continuous at its singularities, we do not establish any facts about $RF$ for $b  \in \left\{\left(\frac{n+1}{n}\right)^2\right\}$. We do not have sufficient computing power to directly apply the bounds and methods of Section \ref{noclassesin89} to these $b$-values.

The graph of $RF(b)$ on a subset of the interval in question appears in Figure \ref{rffigure:bonetotwo}. The class determining the $RF$-value changes when passing through $b\in\left\{\frac{(n+1)^2}{n^2}\right\}$. In the figure we depict this change by shading portions of the graph of $\mu(S_k)$ or $\mu(T_k)$ in gray. Importantly, Jin and Lee's result in \cite{JL} only involved the classes $T_k$ because $\frac{m+1}{m}\in\left[\frac{(n+2)^2}{(n+1)^2},\frac{(n+1)^2}{n^2}\right]$ only for even $n$, hence only the $T_k$ classes were needed to conclude that result.

%In particular $RF(2)$ does not follow the formula above, so already the function $RF:[1,\infty)\to\R$ is not continuous at $b=2$.

% In addition, the result below establishes that $RF(b)$ also fails to be continuous at $b=1$. We show that for the sequence $b_n = \frac{n+1}{n}$, we have $\lim_{n\to\infty} RF(b_n)\neq RF(1)$.
% \begin{thm}\label{theoremB}
% For the sequence $b_n=\frac{n+1}{n}$, with $n\geq 5$,
% \begin{equation}
% RF(b_n) = \frac{2(n+1)}{n}\left(\frac{8n^2+8n+1}{2(2n+1)(n+1)}\right)^2.\label{bnRFvalue}%error corrected here, should be squared
% \end{equation}

% \end{thm}
% When $n\to\infty$, (\ref{bnRFvalue}) approaches 8, whereas  \cite{FrMu} shows that $RF(1) = 7\frac{1}{32}$.

\subsection{Outline of Paper}
In Section 2, we provide some background results needed in later sections, and prove Theorem \ref{mainThm} assuming the contents of the subsequent sections. Then, in Section 3 we describe the exceptional classes which determine the $RF$-value for $b\not\in\left\{\left(\frac{n+1}{n}\right)^2\right\}$, and establish that they are indeed obstructive. In Section 4, we outline the computer-assisted search establishing that no more obstructive classes can appear when $a\in (8,9]$. The remainder of the paper is devoted to applying the reduction method for $a$ in the interval $[9,\infty)$, to establish that there are no obstructions to embeddings for such $a$ other than the volume constraint. This argument splits into many cases, according to the ordering of the terms.

\subsection{Acknowledgments}
The first author would like to thank Washington and Lee University for their support via the Lenfest Grant. The second author thanks Dan Cristofaro-Gardiner for suggesting this problem and providing support and encouragement along the way.

This work used the NCSA Delta CPU at the University of Illinois Urbana-Champaign through allocation MTH230018 from the Advanced Cyberinfrastructure Coordination Ecosystem: Services \& Support (ACCESS) program, which is supported by National Science Foundation grants \#2138259, \#2138286, \#2138307, \#2137603, and \#2138296.

\section{Preliminaries: Two Methods for Finding Symplectic Embeddings}
Here, we review the methods we use for detecting symplectic embeddings. This is not an exhaustive list, and more detailed expositions are in \cite{CFS,McSc}, so we review only what we use. The following definition is central to both methods. Fix $b\geq 1$. Since the function $c_b(a)$ is continuous in $a$, it suffices to compute it for $a \geq 1$ rational.
\begin{defi}\label{weightexpansion}
\textit{The \textbf{weight expansion} $\w(a)$ for $a\in\mathbb{Q}$ is the finite decreasing sequence
\be
\w(a) = (1^{\times \ell_0}, w_1^{\times \ell_1},...,w_n^{\times \ell_n})\nonumber,
\ee
where $w_1 = a - \ell_0 < 1, w_2 = 1 - \ell_1 w_1 < w_1$, and so on.
}
\end{defi}
Mcduff in \cite[Thm. 1.1]{McD}shows that the embedding problems below are equivalent, where $k\in\N$ and $\mathring{B}(\mu)$ is an open ball with radius $\mu$.
\be
E(1,k)\hookrightarrow \mathring{B}(\mu) \iff \coprod_{i=1}^k B(1)\hookrightarrow \mathring{B}(\mu)\nonumber
\ee
In addition, there is a canonical decomposition of the ellipsoid $E(1,a)$ into a finite disjoint union of balls $\coprod_i\coprod_{j=1}^{\ell_i} (w_i)$ where the weights $w_i$ come from the weight expansion of $a$.
Based on this, Frenkel and M\"uller in \cite{FM} to establish the following.
\begin{thm}[{\cite{FM}}]
\textit{Let $a,b\in\Q$ be positive, and let $w_i$ be the terms of the weight expansion $\w(a)$. There exists a symplectic embedding $E(1,a)\hookrightarrow P(\lambda, \lambda b)$ if and only if there is a symplectic embedding}
\be \left(\coprod_i\coprod_{j=1}^{\ell_i} B(w_i)\right)\coprod B(\lambda)\coprod B(\lambda b)\hookrightarrow B(\lambda+\lambda b)
\nonumber
\ee
\end{thm}
This reduces the polydisc problem to a ball-packing problem of embedding balls of capacity $e_i$ into a ball of capacity $\mu$:
\begin{equation}
\coprod_i B(e_i)\hookrightarrow B(\mu).\label{ballembedding}
\end{equation}
%This relates to holomorphic curves via the correspondence between symplectic embeddings of balls and blow-ups. Briefly, if one can embed a ball $B^{2n}(V+\epsilon)$ of volume $V+\epsilon$ into a symplectic manifold $M^{2n}$, then there is a symplectic manifold $\tilde{M}$ obtained from the union of $M\setminus B^{2n}(V+\epsilon)$ and a neighborhood of the zero-section in $\mathcal{O}_{\mathbb{P}^{n-1}}(-1)$. This neighborhood has total volume $\epsilon$ where the zero-section has symplectic volume $V$ \cite[Ch. 7]{McSa}. When $n=2$, this corresponds to an embedded pseudoholomorphic curve, called the exceptional sphere. Specific to this four-dimensional case, let $X_t$ denote the blow-up of $\mathbb{C}P^2$ in $n$ points.\newline\indent
The purpose of introducing these constructions is that a collection of balls can be embedded symplectically precisely when the associated multiple blow-up of $\mathbb{C}P^2$ carries a symplectic form in a certain cohomology class. Let $X_t$ denote the blow-up of $\mathbb{C}P^n$ in $n$ points. We denote by $\overline{\mathcal{C}}_K(X_t)$ the set of cohomology classes represented by symplectic forms for which the anticanonical class is $K=-3L+\sum_i E_i$. Here $L$ is Poincar\'e dual to a line in $\mathbb{C}P^2$ and each $E_i$ is dual to the $i$th exceptional sphere (a symplecti-
cally embedded $S^2$ of self-intersection -1). By \cite{McPo}, the embedding (\ref{ballembedding}) exists when the following cohomology class is in the symplectic cone:
\be
\mu L - \sum_i e_i E_i \in \overline{\mathcal{C}}_K(X_t).\newline\indent
\ee
The above fact gives a sufficient criterion for a class to lie in $\overline{\mathcal{C}}_K(X_t)$. If there is a symplectic form in a given class, it must have non-negative intersection with certain holomorphic curves. By \cite{LiLiu}, this is also sufficient: if $\mathcal{E}_K(X_t):=\{e\in H_2(X_t)\,|\,  \langle e, e\rangle = -1,\,\,\langle K, e\rangle = -1\}$ is the set of exceptional classes, then we may characterize the symplectic cone referenced above as
% -KE = 1 right? so KE = -1.
\begin{equation}
\overline{\mathcal{C}}_K(X_t)=\{\omega\in H^2(X_t)\,|\, \omega^2\geq 0,\, \langle \omega, e\rangle\geq 0\,\, \forall e\in \mathcal{E}_K(X_t)\}. \label{symplecticcone}
\end{equation}
In summary, an embedding of the form (\ref{ballembedding}) exists when the associated weight expansion represents a cohomology class in the symplectic cone $\overline{\mathcal{C}}_K(X_t)$, and the equality (\ref{symplecticcone}) characterizes elements of the symplectic cone as those which pair non-negatively with the exceptional classes $\mathcal{E}_k(X_t)$.

A consequence of this is the following useful alternate characterization of the function $c_b(a)$ in terms of these classes
.

\begin{cor}\label{cbasup}
\textit{\cite[Corollary 1.3]{Ush}For any $b\geq 1$ and any $a\in\mathbb{Q}$ whose weight sequence has length $k-1$ we have}
\be
c_b(a) = \max\left\{\sqrt{\frac{a}{2b}}, \sup_{C\in \mathcal{E}_K} \mu_b(C)(a)\right\}.
\ee
\end{cor}

We now require a condition for determining when homology classes are exceptional. 
\begin{defi}
With respect to the basis of $H^2(X_t;\R)\simeq \R^{n+1}$ given above, a \emph{Cremona transform} is the map given by
\be\label{cremonatransform}
(d; m_1,...m_n)\mapsto (2d-m_1-m_2-m_3;d-m_2-m_3, d-m_1-m_3, d-m_1-m_2, m_4,...,m_n)\nonumber
\ee
\end{defi}
We will use Cremona transforms in two different contexts below. First, we can state the condition for an integral homology class to lie in $\mathcal{E}_K$, which is proven in \cite{McSc} based on work of \cite{LiLi, LiLiu}.
\begin{thm}\label{exceptionalcriterion}
\textit{A class $(d;m_1,...,m_n)\in H_2(X_t;\Z)$ is in $\mathcal{E}_K(X_t)$ if and only if its entries satisfy the Diophantine equations
\begin{eqnarray}
& & 3d-1 = \sum_i m_i \label{diophantine} \newline\\
& & d^2+1 = \sum_i m_i^2 \nonumber
\end{eqnarray}
and $(d;m_1,...,m_n)$ reduces to $(0;-1,0, \ldots, 0)$ after a sequence of Cremona transforms.}
\end{thm}

For the problem of embedding ellipsoids into polydiscs, the natural compactification of $P(a,b)$ adds a single point to each disc, yielding $S^2\times S^2$. There is a diffeomorphism $\psi$ from this manifold to the 2-fold blow-up of $\mathbb{C}P^2$, so the $n$-fold blow-up of $S^2\times S^2$ (denoted $Y_n$) can be identified with $X_{n+1}$. The induced isomorphism on homology $\psi_*:H_2(Y_n)\to H_2(X_{n+1})$ is given by
\begin{equation}
\psi_*(\langle d,e;m_1,m_2,m_3,...,m_n\rangle)= (d+e-m_1; d-m_1, e-m_1, m_2,...,m_n)\label{coordchange}
\end{equation}

Using this isomorphism we may translate the Diophantine equations (\ref{diophantine}) into the following conditions.
\begin{eqnarray}
\sum_i m_i & = & 2(d+e) -1 \label{dioph1}\newline\\
\sum_i m_i^2 & = & 2de+1 \label{dioph2}
\end{eqnarray}

\subsection{Obstructive Classes}\label{method1}
We can now describe the first of the two methods we employ to find symplectic embeddings. The following statement is \cite[Method 1']{CFS}.

\begin{thm}\label{classcondition}
\textit{An embedding (\ref{embeddingproblem}) exists if and only if $\lambda\geq \sqrt{\frac{a}{2b}}$ and
\begin{equation}
\mu_b(d,e;m)(a):= \frac{\sum_i m_i \cdot w_i(a)}{d+be} \leq \lambda\label{obstr}
\end{equation}
for every $(d,e;m_1,...,m_n)\in H_2(Y_n;\Z)$ which satisfies equations (\ref{dioph1}, \ref{dioph2}) and reduces to $(0,1;1)$ after some sequence of Cremona transforms.}
\end{thm}

\begin{defi}
\textit{We say that a class $A := (d,e;\textbf{m})\in \mathcal{E}_K$ is \textbf{obstructive at $a>0$} if $\mu_b(d,e;\textbf{m}){d+be}$
is larger than the volume constraint $\sqrt{\frac{a}{2b}}$.}
\end{defi}
To use this method, then, we must find all obstructive classes. In practice, this method proceeds first by finding a discrete subset of $a$-values which have a special relationship to the obstructive classes, described as follows.

It can be proven as in Corollary 2.1.4 of \cite{McSc} by McDuff and Schlenk that if the graph of $c_b(a)$ does not follow the volume constraint, then locally it must be given by the obstruction function of some class $(d,e;\m)$. Restricting to the interval where this class determines the graph, Cristofaro-Gardiner, Frenkel, and Schlenk show in \cite[Lemma 3.11]{CFS} that for the problem (\ref{embeddingproblem})) there is a particular $a$-value on the interval whose weight expansion coincides with the number of positive entries of $\m$, the tail of the obstructive class.

\begin{lem}\label{lengthbound}
\cite[Lemma 3.11]{CFS}\textit{Let $C = (d,e;\textbf{m})$ be an exceptional class, and let $I$ be the maximal nonempty open interval on which $\mu_{a,b}(C)>\sqrt{\frac{a}{2b}}$ for all $a\in I$. Then, there is a unique element $a_0\in I$ such that the length of the weight expansion for $a_0$, denoted $\ell(a_0)$, satisfies $\ell(a_0) = \ell(\textbf{m})$, and moreover $\ell(a)>\ell(\textbf{m})$ if $a\in I\setminus\{a_0\}$.}
\end{lem}

\begin{defi}\label{centerDef}
\textit{If $C = (d, e; \mathbf{m})$ is an exceptional class, and $a_0 \in I$ is as in Lemma \ref{lengthbound}, we refer to $a_0$ as the \textbf{center} of the obstructive class $C.$}\end{defi}

The following result, \cite[Lemma 2.1.7]{McSc} constrains the multiplicities of entries for candidate obstructive classes.
\begin{lem}\label{oneblock}
\textit{Assume that $(d,e;\m)$ is an exceptional class such that $\mu(d,e;\m)>\sqrt{\frac{a}{2b}}$ for some $a$. Let $J:=\{k,...,k+s-1\}$ be a block of $s$ consecutive integers for which $w(a_i), i\in J$ is constant. Then $(m_k,\dots, m_{k+s-1})$ is of the form}
\bea
(m,...,m)\newline\nonumber\\
(m,...,m,m-1),\mbox{or}\newline\nonumber\\
(m+1,m,...,m).\nonumber
\eea
Moreover, there is at most one block of length $s\geq 2$ on which the $m_i$ are not all equal, and $\sum_{i=k}^{k+s-1}\epsilon_i^2 \geq \frac{s-1}{s}$.
\end{lem}
One more lemma below further restricts the shape of candidate obstructive classes by limiting their possible centers. In order to state this lemma we will need to define the following quantities. Let $\ell_0$ be as in Definition \ref{weightexpansion} and subsequently let $\ell_i$ denote the lengths of subsequent blocks. When a class $(d,e;\w(a))$ is obstructive, we have a vector of error terms $\epsilon$ defined as
$\m = \frac{d+be}{\sqrt{2a}}\w(a)+\boldsymbol{\epsilon}$.
Each contribution to this error, thought of as the difference between $\m$ and $\frac{d+be}{\sqrt{2a}}\w(a)$, will be important, so we define $v_i = \frac{d+be}{\sqrt{2a}}w_i$ for $i=0,...,M$ where $M$ denotes the length of the weight expansion of $a$.\newline\indent
Also let
\be
\sigma = \sum_{\ell_0 + 1}^M \epsilon_i^2. \nonumber
\ee
% denote the ``residual error," the contribution to the error vector coming from non-integer terms of $\w(a)$.
Also, define a related quantity
\be
\sigma' = \sum_{\ell_0 + 1}^{M-\ell_N} \epsilon_i.\nonumber
\ee
% Here $\ell_N$ is the length of the last block. This quantity ignores the contribution to the error from the smallest part of the weight expansion $\w(a)$.  \newline\indent
With these terms established, the following result \cite[Proposition 3.1(iii)]{CFS} will be used repeatedly.
\begin{thm}\label{obstructiveconditionsh}
Suppose that $(d,e;\m)$ satisfies the Diophantine equations (\ref{dioph1}), (\ref{dioph2}).
\begin{enumerate}
    \item $\mu_b(d,e;m) >\sqrt{\frac{a}{2b}}$ iff the error vector satisfies $\langle \epsilon, \textbf{w}(a)\rangle >0$
    \item If $\mu_b(d,e;m) >\sqrt{\frac{a}{2b}}$ then $d = be + h$ where $|h| < \sqrt{2b}$, and $||\epsilon||^2 < 1-\frac{h^2}{2b}$. \label{errorbound}
\end{enumerate}
\end{thm}

The results below are from (3) and (4) of \cite[Lemma 2.8]{JL}. The inequality (\ref{ebound}) is a bound on the magnitude of the entries of candidate obstructive classes, as an upper bound on $e$ limits the right-hand sides of (\ref{dioph1}) and (\ref{dioph2}).

\begin{lem}\label{errorestimates}%Analogue of CFS Lemma 3.7
\textit{Let $(d,e;\textbf{m})$ be an exceptional class such that $\ell(a) = \ell(m)$ and $\mu_b(d,e;\m) > \sqrt{\frac{a}{2b}}$ for some $a\in (8,9)$, and $b\in[1,2]$.
Set $v_M := \frac{(d+be)\sqrt{2b}}{q(b+1)\sqrt{a}}$ where $q$ is the last denominator in the weight expansion $\w(a)$. Then}
\begin{enumerate}
%   \item $\left|\sum_i \epsilon_i\right| \leq \sqrt{\sigma L}$
%   \item If $v_M<1$, then $\left| \sum_i \epsilon_i\right|\leq \sqrt{\sigma'L}$,
   \item \textit{If $v_M\leq \frac{1}{2}$, then $v_M> \frac{1}{3}$ and $\sigma'\leq \frac{1}{2}$. If $v_M\leq \frac{2}{3}$, then $\sigma'\leq \frac{7}{9}$.}\label{vmsigmabounds}
   \item \textit{For $\delta=$ $a+1-2\frac{b+1}{\sqrt{2b}}\sqrt{a}-\frac{1}{q}$, we have},
\begin{equation}
2be+h \leq \frac{\sqrt{2ba}}{\delta}\left(\sqrt{\sigma q} - \left(1-h\left(1-\frac{1}{b}\right)\right)\right) \leq \frac{\sqrt{2ba}}{\delta}\left(\frac{\sigma}{\delta v_M}-\left(1-h\left(1-\frac{1}{b}\right)\right)\right). \label{ebound}
\end{equation}

% \item  \label{qbound}
%   \begin{equation}\label{qboundformula}
%   \sqrt{q + \floor{a} + 1} \geq 1 + \delta v_M q
%   \end{equation}
\end{enumerate}
\end{lem}

% The inequality (\ref{qbound}) restricts the denominators $q$ which may appear as centers $\frac{p}{q}$ of obstructive classes, which in turn restricts the possible shapes of such classes.

% \begin{lem}\label{eboundlemma}
% For $\langle d,e;\m\rangle$ satifying the hypotheses of Lemma \ref{errorestimates} at $a\in(8,9)$,
% \begin{equation}
% 2be+h \leq \frac{\sqrt{2ba}}{\delta}\left(\sqrt{\sigma q} - \left(1-h\left(1-\frac{1}{b}\right)\right)\right) \leq \frac{\sqrt{2ba}}{\delta}\left(\frac{\sigma}{\delta v_M}-\left(1-h\left(1-\frac{1}{b}\right)\right)\right). \label{ebound}
% \end{equation}
% \end{lem}

\subsection{Reduction at a Point}\label{method2}
Though the first method is a necessary and sufficient condition for the embedding of an ellipsoid into a polydisc, in principle one might need to check more exceptional classes than is computationally feasible. The following method straightforwardly determines whether an embedding exists but the complexity of the process depends strongly on the value of $a$.\newline\indent

\begin{defi}
\textit{The \textbf{defect} $\delta$ of an ordered vector $(d;m_1,...,m_n)$ is the sum $d-m_1-m_2-m_3$. A vector is \textbf{reduced} if $\delta\geq 0$.}
\end{defi}
Buse-Pinsonnault and Karshon-Kessler establish following in \cite{BuPi, KaKe}.
\begin{thm}\label{reductionThm}
\textit{An embedding $E(1,a)\hookrightarrow P(\lambda, \lambda b)$ exists if there exists a finite sequence of Cremona moves that transforms the ordered vector
\begin{equation}
((b+1)\lambda; b\lambda, \lambda, \w(a))\label{reductionvector}
\end{equation}
to an ordered vector with non-negative entries and defect $\delta\geq 0$.}
\end{thm}
We will apply this repeatedly in the proof of Theorem 1.1.

Notice the requirement that the Cremona moves be performed on ordered vectors. The following fact \cite[Prop. 2.2]{CFS} allows us to avoid completely determining the order of all entries in a given vector.

\begin{prop}
\textit{Let $\alpha=(\mu;a_1,...,a_n)$ be a vector with $\mu\geq 0$ and $\langle\alpha,\alpha\rangle\geq 0$ and assume that there is a sequence $\alpha = \alpha_0,...,\alpha_m$ of vectors such that $\alpha_{j+1}$ is obtained from $\alpha_j$ by a sequence of Cremona moves. If $\alpha_m = (\hat{\mu};\hat{\alpha}^1,...,\hat{\alpha}^n)$ is reduced and $\hat{\alpha}_1,...,\hat{\alpha}_n\geq 0$, then $\alpha\in\overline{\mathcal{C}}_K(X_t)$.}
\end{prop}

Thus to find embeddings we need only prove enough about the ordering to ensure that the defect at a certain step is non-negative, provided that we have non-negativity of the terms in each vector. So, only knowledge of the ordering for the first three $a_i$-terms is strictly necessary once we know these must be the largest three. We will usually verify the ordering at each step, except when doing so results in distinguishing too many cases.

\subsection{Proof of Theorem \ref{mainThm}}

At this point we are ready to write a proof of Theorem \ref{mainThm}, using the results of sections which follow this one. We restate it below for convenience.

\begin{customthm}{1.1}
\textit{
Let $b \in \left(1, 2 \right)\setminus\left\{\frac{(n+1)^2}{n^2}: n \in \mathbb N\right\}.$ Then
\begin{equation}
RF(b) =
\begin{cases}
    2b(\mu_b(S_{k+1})(8))^2 & \mbox{if $b \in I_{2k+1}$} \\
    2b(\mu_b(T_k)(8))^2 & \mbox{if $b \in I_{2k}$} \\
   \end{cases}.
\end{equation}}

\end{customthm}

\begin{proof}
Fix $b\in I_{2k}$; the argument for $b\in I_{2k+1}$ is similar. Define
\be
X := \{A \in \mathbb R: c_b(a) = \sqrt{\frac{a}{2b}}\text{ for all } a \ge A\},
\ee
To see that $RF(b) = 2b(\mu_b(T_k)(8))^2 =: \rho$ is a lower bound of $X$, we will show that for $a\in[8, \rho)$ we have $c_b(a) > \sqrt{\frac{a}{2b}}$. Note that indeed $\rho > 8$ by Lemma \ref{sktksupports}. Also by
Theorem \ref{cbaateight}, we know that $c_b(8) = \mu_b(T_k)(8) > \frac{2}{\sqrt{b}}$. By the definition of $c_b(a)$ we have, for all $a$ and any exceptional class $E$, that $c_b(a) \geq \mu_b(E)(a)$. In particular, for $a \geq 8$ we have $c_b(a) \geq \mu_b(T_k)(a)$. We claim, and will prove next, that for all $a \geq 8$ we have $\mu_b(T_k)(a) = \mu_b(T_k)(8)$, and therefore $c_b(a) \geq \mu_b(T_k)(8)$. Finally, we note that by definition of $\rho$ the curve $\sqrt{\frac{a}{2b}}$ intersects the horizontal line of height $\mu_b(T_k)(8)$ exactly at $a = \rho$, and thus for $a <\rho$ we have $\mu_b(T_k)(8) > \sqrt{\frac{a}{2b}}$ , and therefore $c_b(a) > \sqrt{\frac{a}{2b}}$ for
$a \in [8, \rho)$. We now prove the earlier claim that for $a \geq 8$ we have $\mu_b(T_k)(a) = \mu_b(T_k)(8)$:
this holds because $\ell(T_k)=8$, so further terms in $\w(a)$ for $a> 8$ do not contribute.

To see that in fact $\rho\in X$, first note that by Theorem \ref{agreaterthan9}, we have that $c_b(a)=\sqrt{\frac{a}{2b}}$ for $a\geq 9$. For $a_0\in (\rho,9)$, if $c_b(a_0)>\sqrt{\frac{a_0}{2b}}$, then by continuity of $c_b(a)$ in $a$ there exists a maximal open interval $U$ containing $a_0$ on which this inequality holds. This implies that there exists an $E\in \mathcal{E}_K$ which is obstructive at $a_0$, and whose center is some $a'$.

By Proposition \ref{noclassesin89}, $a'\leq 8$, so $8\in U$. Hence $E$ is obstructive at $a=8$, so by Lemma \ref{onlycenteredat8}, $E$ must be one of the $S_j$ or $T_j$ for some $j\geq 1$, or $E=\langle 3,1;1^{\times 7}\rangle$. But by Lemmas \ref{sktksupports} and \ref{obstructiveClassLemma}, and the assumption that $b\in I_{2k}$, $E$ must be of the form $T_k$, and moreover the support of $\mu_b(T_k)(a)$ lies below $\rho$. So $c_b(a) = \sqrt{\frac{a}{2b}}$ for $a\geq \rho$ as needed.

\end{proof}

\section{The classes \texorpdfstring{$S_k, T_k$}{Sk, Tk} determine \texorpdfstring{$c_b(8)$}{cb8} for \texorpdfstring{$b\not\in\left(\frac{(2k+2)^2}{(2k+1)^2},\frac{(2k+1)^2}{(2k)^2}\right)$}{b not a ratio of consecutive squares}}\label{sktkprops}

In this section we show that the value of the ellipsoid embedding function $c_b(8)$ is determined by the  classes (\ref{skclasses}), (\ref{tkclasses}).

\begin{lem}\label{onlycenteredat8}
If $C = \langle d, e; \m \rangle$ is obstructive at $a = 8,$ then either $C$ is centered at $a = 8$ or $C = \langle 3,1;1^{\times 7}\rangle$.

\end{lem}
\begin{proof}
Let $C = \langle d, e; \mathbf{m}\rangle$ be an obstructive class at $a = 8.$ That is, suppose $C \in \mathcal{E}_K$ and $\mu_b(C)(8) > \frac{2}{\sqrt{b}}$ on a maximal open interval $J.$ By Lemma \ref{lengthbound}, there exists a unique $a_0 \in J$ such that $\ell(a_0) = \ell(\m),$ and if $a \in J\setminus\{a_0\},$ then $\ell(a) > \ell(\m).$

Either $a_0 = 8$ or $a_0 \ne 8.$ If $a_0 = 8,$ then $\ell(\m) = \ell(8) = 8.$ Otherwise, if $a_0 \ne 8,$ then by Lemma \ref{lengthbound}, $8 = \ell(8) > \ell(\m)$ and hence $\ell(\m) \le 7.$ Now, Frenkel and M\"uller show in \cite[Lemma 4.3]{FM} there are only ten classes in $\mathcal{E}_K$ such that $\ell(\m) \le 7$.

Now if $\ell(\m)\le 7$, then at least one $m_i=0$. As $C$ is obstructive at $a=8$, the first eight terms of $\m$ must be one of the three possibilities listed in Lemma \ref{oneblock}: in the notation of that lemma, either $m=0, m+1=0$, or $m-1=0$. Thus, out of the ten classes with length at most seven, only $\langle 3,1;1^{\times 7}\rangle,\langle 0,0;-1\rangle, \langle 1,0;1\rangle$ satisfy this condition, and the classes containing 0 are not obstructive.
\end{proof}

Now we show that the only possible obstructive classes centered at $a = 8$ are of the form $S_k$ or $T_k$ for some $k \ge 1.$ Specifically, we seek to prove the following result:

\begin{lem}\label{obstructiveClassLemma} \textit{For $1 < b \le 2,$ the only exceptional classes $\langle d, e ; \mathbf m \rangle $ with $\mu_b(d, e; \mathbf m)(8) > \frac{2}{\sqrt{b}}$ are $S_k$ and $T_k$ for $k \ge 1.$} \end{lem}

First, we note that if $C = \langle d, e; \m \rangle$ is obstructive at $a = 8,$ Lemma \ref{KeyLemma} (which is simply a restatement of Lemma \ref{oneblock}),determines the possible shapes of $\m:$

\begin{lem}\label{KeyLemma} \textit{If $\langle d, e; \mathbf{m}\rangle$ is obstructive with center at $a = 8,$ then $\mathbf{m} = (m^{\times 8}), (m^{\times 7}, m-1)$ or $(m+1, m^{\times 7}).$}\end{lem}

In addition, recall that if $\langle d, e; \m\rangle$ is obstructive, where $\m = (m_1, \ldots, m_k),$ then

\begin{equation}\tag{\ref{dioph1}}
    2de + 1 = \sum_{i=1}^k m_i^2,
\end{equation}
and
\begin{equation}\tag{\ref{dioph2}}
2(d+e) - 1 = \sum_{i=1}^k m_i.
\end{equation}

Finally, a basic number-theoretic result, whose proof is straightforward and we thus omit, will be useful in finalizing the classification:
\begin{lem}\label{square} \textit{Let $x \in \mathbb Z.$ Then $8|(x^2 - 1)$ if and only if $x = \pm 1 \pmod 4.$}\end{lem}
%\begin{proof} If $x = \pm 1+ 4k$ for $k \in \mathbb Z,$ then \be x^2 = 1 \pm 8k + 16k^2,\ee so $8|(x^2 - 1).$ Conversely, if $x^2 = 1 \pmod 8,$ then $x^2 = 1 \pmod 4,$ so $x \ne 0, 2 \pmod 4,$ therefore $x = \pm 1 \mod 4.$ \end{proof}

With these ideas in mind, we are ready to prove Lemma \ref{obstructiveClassLemma}.

\subsection{Proof of Lemma \ref{obstructiveClassLemma}}
In the definition (\ref{obstr}) of $\mu_b(d,e;m)$, with $b\geq 1$ and $d<e$, $\mu_b(e,d;m)>\mu_b(d,e;m)$, so we may assume $d \geq e.$ Suppose $\langle d, e; \mathbf{m} \rangle$ is obstructive at $a = 8.$ If $\mathbf{m} = (m^{\times 8}),$ then $2de + 1 = 8m^2$ by Equation \ref{dioph1}, which is a contradiction because $d, e$ and $m$ are all integers. Therefore, the only two cases we need to consider are $\mathbf{m} = (m^{\times 7}, m - 1)$ and $\mathbf{m} = (m+1, m^{\times 7}).$

\subsubsection{Case I: \texorpdfstring{$\mathbf{m} = (m^{\times 7}, m - 1)$}{7 ms and one m-1}}

In this case, we have:
\begin{equation}\label{C1E1}
2de+1 = 8m^2 - 2m + 1,
\end{equation}
and
\begin{equation}\label{C1E2}
2(d + e) - 1 = 8m - 1.
\end{equation}

By Equation \ref{C1E2}, $m = (d + e)/4.$ Substituting $m$ into Equation \ref{C1E1}, we see that

\begin{equation}\label{deQ}
de = \frac{(d+e)^2}{4} - \frac{d+e}{4}.
\end{equation}

Let $h = d + e.$ Then $de = (d+e - e)e = he - e^2.$ Upon clearing denominators in \ref{deQ} and viewing the result as quadratic in $h,$ we have
\begin{equation}\label{2H}
2h = 1 + 4e \pm \sqrt{8e + 1}.
\end{equation}

Since $h \in \mathbb Z,$ it follows that $8e+1$ is a perfect square. So there exists $x \in \mathbb Z$ such that $x^2 = 8e + 1.$ Hence $x^2 = 1 \pmod 8,$ so $x = \pm 1 + 4k$ for some $k \in \mathbb Z$ by Lemma \ref{square}.
%In particular, if $x^2 = 8e + 1$ for some $x \in \mathbb Z,$ then, by Lemma \ref{square}, $x = \pm 1 + 4n$ for some $n \in \mathbb Z.$
Thus, either $e = 2k^2  + k$ or $e = 2k^2 - k.$

%On the other hand, \be 8e+1 = (\pm 1 + 4n)^2 = (1 + 4n)^2,\ee so \begin{equation}\label{new2H} 2h = 8n^2 + 4n + 1 \pm (1 + 4n).\end{equation}$}

If $e = 2k^2 - k,$ then $8e + 1 = (4k - 1)^2,$ and we have $2h = 8k^2 - 4k + 1 \pm (4k - 1).$ So either $h = 4k^2$ or $h = 4k^2 - 4k + 1.$ The latter case cannot happen as $4|h,$ so we must have $h = 4k^2.$ Therefore, if $e = 2k^2 - k$ for some $k \in \mathbb Z,$ then \be d = h - e = 2k^2 + k.\ee

%Thus, the only possible class arising from this case is \begin{equation}\label{goodsnClasses}
%S_n:= \langle 2n^2 + n, 2n^2 - n; (n^2)^{\times 7}, (n^2 - 1)\rangle \end{equation}

If $e = 2k^2 + k,$ then $8e+1 = (4k+1)^2$, and like before we have $h = 4k^2$ so that $d = 2k^2 - k.$ However, since we are assuming that $d \ge e,$ this solution is excluded.

In summary, if $\langle d, e; \mathbf{m} \rangle $ is obstructive at $a = 8,$ where $\mathbf{m} = \langle m^{\times 7}, m - 1 \rangle,$ then \be \langle d, e; \mathbf{m} \rangle = \langle 2k^2 + k, 2k^2 - k; (k^2)^{\times 7}, (k^2 - 1)\rangle = S_k.\ee Since $d \ge e,$ we have $2k^2 + k \ge 2k^2 - k$ which implies $k \ge 0.$

\subsubsection{Case II: \texorpdfstring{$\mathbf{m} = (m+1, m^{\times 7})$}{m plus 1 and 7 ms}}

The argument in this case is very similar to the previous case. By Equations \ref{dioph1} and \ref{dioph2}, we have:

\begin{equation}\label{C2E1}
de = 4m^2 + m,
\end{equation}
and
\begin{equation}\label{C2E2}
d + e  = 4m + 1.
\end{equation}

By Equation \ref{C2E1}, $de = m(4m + 1) = m(d + e).$ Solving for $m$ in Equation \ref{C2E2} and substituting the result into the latter equation, we see that
\begin{equation}\label{C2E3}
4de = (d + e - 1)(d + e).
\end{equation}

Let $h = d + e.$ Rewrite Equation \ref{C2E3} in terms of $h$ to obtain
\begin{equation}\label{C2E4}
h^2 - (1 + 4e)h + 4e^2 = 0.
\end{equation}

Therefore, \begin{equation}\label{C2E5} 2h = 1 + 4e \pm \sqrt{1 + 8e}, \end{equation}

which is exactly Equation \ref{2H}.

%Since $h \in \mathbb N,$ it follows that $1 + 8e$ is a perfect square. That is, $8e = z^2 - 1$ for some $z \in \mathbb Z,$ so $e \in\{k(2k + 1): k \in \mathbb Z\} \cup \{k(2k - 1): k \in \mathbb Z\}$ by Lemma \ref{square}.S

Therefore, we have $e = 2k^2 + k,$ or $e = 2k^2 - k.$ We conduct our analysis in a way completely analogous to what we did in Case I. Suppose $e = 2k^2 + k.$ Then $h = 4k^2$ or $h = 4k^2 + 4k + 1.$
In this scenario, the first case is impossible because $h = d + e = 1 \pmod 4.$ Therefore, \be d + e = 4k^2 + 4k + 1,\ee and thus \begin{equation}\label{dCaseII}
d = 4k^2 + 4k + 1 - e = 2k^2 + 3k + 1 = (2k + 1)(k+1).
\end{equation}

Just like above, the case where $e = 2k^2 - k$ is excluded because $d \ge e.$ Now, $m = k^2 + k$ by Equation \ref{C2E2}, so we see that if $\langle d, e; \mathbf{m} \rangle$ is obstructive at $a = 8$ for some $b \ge 1,$ where $\mathbf{m} = (m + 1, m^{\times 7}),$ then we must have \be\langle d, e; \mathbf{m} \rangle = \langle (2k + 1)(k + 1), k(2k + 1); (k^2 + k + 1), (k^2 + k)^{\times 7} \rangle = T_k\ee for some $k \in \mathbb Z.$ As before, since $d \ge e,$ we have $2k^2 + 3k + 1 \ge 2k^2 + k$ which implies $k \ge 0.$ \qed

\begin{lem}\label{sktksupports}
Let $b \in (1, 2).$ We have the following:
\begin{enumerate}
    \item $\mu_b(T_k)(8) > \frac{2}{\sqrt{b}}$ if and only if $b \in I_{2k}.$
    \item $\mu_b(S_{k+1})(8) > \frac{2}{\sqrt{b}}$ if and only if $b \in I_{2k+1}.$
    \item If $b \in I_{2k},$ then $\mu_b(T_k)(8) \ge \mu_b(C)(8)$ for every $C \in \mathcal{E}_K.$
    \item If $b \in I_{2k+1},$ then $\mu_b(S_{k+1})(8) \ge \mu_b(C)(8)$ for every $C \in \mathcal{E}_K.$
\end{enumerate}
\end{lem}

\begin{proof}
To see (1), consider for each $k \ge 1$ the function $f_{T_k}$ defined below:
\begin{eqnarray*}
f_{T_k}(b) & := & \mu_b(T_k)(8)-\frac{2}{\sqrt{b}} = \frac{8k^2+8k+1}{(2k+1)(k+1) + b\cdot k(2k+1)} - \frac{2}{\sqrt{b}} \\
% f_{S_k}(b) & := & \mu_b(S_k)(8)-\frac{2}{\sqrt{b}} = \frac{8k^2-1}{2k^2+k+ b\cdot (2k^2-k)} - \frac{2}{\sqrt{b}}. \\
\end{eqnarray*}

Observe that $f_{T_k}(b) > 0$ if and only if \begin{equation}\label{sqrtb}
(8k^2 + 8k +1)\sqrt{2b} - \sqrt{8}((2k+1)(k + 1) + b(k(2k+1)) > 0.
\end{equation} Consider the polynomial \be P(t) = (-4k^2-2k)t^2 + (8k^2+8k+1)t + (-4k^2-6k-2).\ee Note that if $t = \sqrt{b},$ we obtain the expression in \ref{sqrtb} (scaled by $1/\sqrt{2}$). The roots of $P$ are $(2k+2)/(2k+1)$ and $(2k+1)/2k.$ Since the leading coefficient of $P$ is negative, and the roots of $P$ are distinct, it follows that $P$ is positive if and only if \be\frac{2k+2}{2k+1} < t < \frac{2k+1}{2k}.\ee Equivalently, since $b>1$ by hypothesis, the expression in \ref{sqrtb} is positive if and only if $\left(\frac{2k+1}{2k+1}\right)^2 < b < \left(\frac{2k+1}{2k}\right)^2,$ or equivalently if and only if $b \in I_{2k}.$

Similarly for (2), the analogous function $f_{S_{k+1}}(b):= \mu_b(S_{k+1})(8)-\frac{2}{\sqrt{b}}>0$ if and only if $b \in I_{2k+1}.$

To see that the obstruction from $T_k$ exceeds that of any other $\{S_i\}\cup \{T_i\}_{i\neq k}$ for these values of $b$, take any $b\in I_{2k}$ and consider any other class $C\in\{S_i\}\cup \{T_i\}_{i\neq k}$:

\begin{eqnarray*}
\mu_b(T_k)(8) - \mu_b(C)(8) & = &  \mu_b(T_k)(8) - \mu_b(C)(8) + \frac{2}{\sqrt{b}} - \frac{2}{\sqrt{b}} \\
& = & f_{T_k}(b) - f_C(b).
\end{eqnarray*}

By the preceding argument, $T_k$ is the unique class such that $f_{T_k}(b)>0$ for $b\in I_{2k}$, so $f_{T_k}(b) - f_C(b) > 0.$ Finally, checking that $\mu_b(\langle 3, 1; 1^{\times 7}\rangle)(8)$ is less than volume on $I_{2k}$ amounts to verifying that \be\frac{2}{\sqrt{b}} - \frac{7}{3+b} > 0,\ee which is straightforward. Consequently, $\mu_b(T_k)(8) \ge \mu_b(C)(8)$ for all $C \in \mathcal{E}_K$ on $I_{2k}.$ The argument for $S_{k+1}$ is identical.
\end{proof}

We now verify that $S_k$ is exceptional for each $k \ge 1.$ The proof is similar to \cite[Lemma 7.4]{JL}, where it is shown that each $T_k$ is exceptional (in that work, $R_n$ is used to denote what we call $T_n$ here.) Recall that by Theorem \ref{exceptionalcriterion}, a class $C$ is exceptional if it satisfies Equations (\ref{dioph1}), (\ref{dioph2}) and, after applying the isomorphism \ref{coordchange}, reduces to $(0;-1)$ after a sequence of Cremona transforms.

\begin{prop}\label{Skexceptional}
If $k > 0,$ then $S_k$ reduces to $(0; -1).$
\end{prop}
\begin{proof}
Applying the isomorphism $\psi_*$ from \eqref{coordchange} to \be S_k =  \langle 2k^2 + k, 2k^2 - k; (k^2)^{\times 7}, (k^2 - 1) \rangle,\ee we obtain
\begin{equation}\label{psiStar}
    \psi_*(S_k) = (3k^2;k^2 + k, (k^2)^{\times 6}, k^2 - 1, k^2 - k),
\end{equation} with defect $\delta = -k.$ Applying a standard Cremona move to \eqref{psiStar}, we obtain \begin{equation}\label{Ck}
C_k:=(3k^2 - k; (k^2)^{\times 5}, k^2 - 1, (k^2 - k)^{\times 3}).
\end{equation}

It suffices to show that $C_k$ reduces to $(0; -1)$ via a finite sequence of Cremona moves. If $A$ and $B$ are ordered vectors, we will write $A \to B$ if there exists a finite sequence of Cremona moves (i.e., computation of defect, application of defect, then reordering if necessary) that carries $A$ to $B.$

We argue by induction on $k.$ If $k = 1,$ a direct computation shows that $C_1 \to (0; -1),$ so let $k > 1$ and assume that $C_{k-1} \to (0;-1).$ We show that after four Cremona moves, $C_k \to C_{k-1},$ so that by the induction hypothesis $C_k \to (0; -1).$ To keep track of the intermediate vectors that will carry us from $C_k$ to $C_{k-1},$ we will enumerate $E_0:=C_k,$ and let $E_i$ be the vector obtained after doing a Cremona move to $E_{i-1}.$

The defect of $E_0 = C_k$ in \eqref{Ck} is $\delta = -k.$ Applying a Cremona move to $C_k$, we have \be E_1 = (3k^2 - 2k; (k^2)^{\times 2}, k^2 - 1, (k^2 - k)^{\times 6}),\ee with defect $\delta = 1 - 2k.$ Again applying a Cremona move, we have \be E_2 = (3k^2-4k+1; (k^2 - k)^{\times 6}, (k^2 - 2k + 1)^{\times 2}, k^2 - 2k),\ee with defect $\delta = 1 - k.$ Therefore, \be E_3 = (3k^2 - 5k + 2; (k^2 - k)^{\times 3}, (k^2 - 2k + 1)^{\times 5}, k^2 - 2k),\ee with defect $\delta = 2 - 2k.$ Finally, with one more Cremona move, we have

\begin{eqnarray*}
E_4 &=& (3k^2 - 7k + 4; (k^2 - 2k + 1)^{\times 5}, k^2 - 2k, (k^2 - 3k + 2)^{\times 3})\\
&=& (3(k-1)^2 - (k-1); ((k-1)^2)^{\times 5}, (k-1)^2 - 1, ((k-1)^2 - (k-1))^{\times 3})\\
&=& C_{k-1}.
\end{eqnarray*}
Therefore, $C_k \to C_{k-1} \to (0;-1),$ where the latter reduction follows from the induction hypothesis. \end{proof}

\begin{thm}\label{cbaateight}
Let $b \in (1, 2)\setminus\left(\frac{(2k+2)^2}{(2k+1)^2},\frac{(2k+1)^2}{(2k)^2}\right)$. Then
\be c_b(8) =
\begin{cases}
\mu_b(T_k)(8) &\text{if } b \in I_{2k},\\
\mu_b(S_{k+1})(8) &\text{if } b \in I_{2k+1.}%\\
%\frac{2}{\sqrt{b}} &\text{otherwise.}
\end{cases}\ee
\end{thm}

\begin{proof}
If $b\in I_{2k}$, then by Lemma \ref{sktksupports} (1) and Proposition \ref{Skexceptional}, $T_k$ is obstructive. Moreover, by Lemma \ref{sktksupports} (3), $\mu_b(T_k)(8)\geq \mu_b(C)(8)$ for any $C\in \mathcal{E}_K$, so $c_b(8) = \mu_b(T_k)(8)$. The proof for $b\in I_{2k+1}$ is similar.
\end{proof}

\section{No obstructive classes for \texorpdfstring{$a\in(8,9]$}{a between 8 and 9}}
In this section we outline the computer-assisted search method used to investigate the interval $(8,9)$.

\begin{prop}\label{noclassesin89}
For $b\in (1,2)$ in the embedding problem (\ref{embeddingproblem}), there are no obstructive classes with center on $(8,9)$.
\end{prop}

\begin{proof}
The argument here is as in \cite[Proposition 7.2]{JL} and Cristofaro-Gardiner-Frenkel-Schlenk \cite[Lemma 3.9]{CFS}. We explain its application briefly here. 

Suppose that $C = \langle d,e;\m\rangle \in \mathcal{E}_K$ is obstructive and has center $a_0\in (8,9)$.  Lemma \ref{lengthbound} requires that $\ell(\m) = \ell(a_0)$, and this length in turn depends on how we may write $a_0=\frac{p}{q}$. The inequality (\ref{ebound}) from Lemma \ref{errorestimates} produces bounds on both the denominator $q$ and the value $e$ appearing in $\langle d,e;\m\rangle$ in as follows: since any $C$ as above must satisfy the inequality \ref{ebound}, recall the sides of this inequality depend on $q,h,b$ and $a\in (8,9)$:
\be
\frac{\sqrt{2ba}}{\delta}\left(\sqrt{\sigma q} - \left(1-h\left(1-\frac{1}{b}\right)\right)\right) \leq \frac{\sqrt{2ba}}{\delta}\left(\frac{\sigma}{\delta v_M}-\left(1-h\left(1-\frac{1}{b}\right)\right)\right).
\ee
Now, for fixed $q,h$ it is straightforward to verify that these are decreasing functions in $a$, so it suffices to consider $a$-values of the form $8\frac{1}{q}$. This reduces to the functions below, which we denote by $u(q,h,b)$ and $l(q,h,b)$:

\bea
l(q,h,b) & = & \frac{\sqrt{2b\cdot 8\frac{1}{q}}}{\delta(q,b)}\left(\sqrt{\sigma q} - \left(1-h\left(1-\frac{1}{b}\right)\right)\right)\\ % =  \frac{\sqrt{2ba}}{\delta(q)}\left(\sqrt{\sigma q} - \left(1-h\left(1-\frac{1}{b}\right)\right)\right)\\
u(q,h,b) & = & \frac{\sqrt{2b\cdot 8\frac{1}{q}}}{\delta(q,b)}\left(\frac{\sigma}{\delta(q,b) v_M}-\left(1-h\left(1-\frac{1}{b}\right)\right)\right).% = \frac{\sqrt{2ba}}{\delta(q)}\left(\frac{\sigma}{\delta v_M}-\left(1-h\left(1-\frac{1}{b}\right)\right)\right)
\eea

As in \cite[Lemma 7.3]{JL}, we have $\frac{\partial f}{\partial q}>0$ and $\frac{\partial g}{\partial q}<0$ for $q>1$, and $\frac{\partial f}{\partial h}, \frac{\partial g}{\partial h}>0$. It is straightforward to verify that for $1<b<2$, and $h^2<2b(1-||\epsilon||^2)$ using (\ref{errorbound}) of Theorem \ref{obstructiveconditionsh}, the largest possible point of intersection of the graphs of $l(q,h,b)$ and $u(q,h,b)$ is approximately $(11.09, 40.56)$. (In \cite[Proposition 7.2]{JL} this argument was used only for $b\in (1,\frac{10}{9}]$ to conclude $q<3$ but its proof applies here to obtain $q<11$ instead.)

With $q=11$ as the largest denominator $q$ in $8\frac{p}{q}$, it is straightforward to verify that there are 41 numbers of the form $8\frac{p}{q}$ with the fraction in lowest terms and $p<q$. This gives rise to 41 possible centers, e.g. $8\frac{1}{2}, 8\frac{1}{3}, 8\frac{2}{3}$, $8\frac{1}{4}$, $8\frac{3}{4}$, and so on.

Next, to use the bound $e\leq 40$, the equation $d=be+h$ and the inequality in Lemma \ref{errorbound} establishes that for each $e$ the range of possible $d$ which might appear in a candidate obstructive class ranges from $e-2$ to $2e+2$. Thus, for each possible $1\leq e \leq 40$ we may list the possible matching $d$ in $\langle d,e;\m\rangle$ satisfying these inequalities. That is, for $e=1$ the possible pairings are $(1,1), (2,1), \dots, (5,1)$. For $e=2$ we have $(1,2)$ to $(6,2)$. Continuing up to $e=40$ in this way we obtain 934 possible $(d,e)$-pairs.

It then remains to verify that no classes with these values of $(d,e)$ are obstructive. The Algorithm \ref{extendalg} below produces all classes $\langle d,e;\m\rangle$ satisfying the conditions of Lemmas \ref{lengthbound} and \ref{oneblock}, from which we may verify that none exceed the volume constraint for $a\in (8,9), b\in (1,2)$.

For an example computation, consider the class $(42, 22; 16, 15^{\times 7}, 3, 3)$. While this satisfies the conditions of Lemma \ref{oneblock} and Lemma \ref{lengthbound} at $a=8\frac{1}{2}$, we can easily calculate that for any $b\in (1,2)$,
\be
\frac{\langle (16, 15^{\times 7}, 3,3), (1^{\times 8}, \frac{1}{2}, \frac{1}{2})\rangle}{42+b\cdot 22} < \frac{\langle (16, 15^{\times 7}, 3,3), (1^{\times 8}, \frac{1}{2}, \frac{1}{2})\rangle}{42+2\cdot 22} < \sqrt{\frac{8\frac{1}{2}}{2\cdot 1}} < \sqrt{\frac{8\frac{1}{2}}{2\cdot b}}.
\ee
\end{proof}

To produce the classes needed to verify our claim, we perform a search along the lines of what appears in \cite[Appendix B]{McSc}. Our search is significantly less economical (to an extreme degree) but it does not depend at all on a choice of $b$. We outline the approach below.

Note that Equation \eqref{dioph1} amounts to the statement that \textbf{m} is a partition of the positive integer $2(d+e)-1$. These are well-understood combinatorial objects, and the problem of efficiently enumerating partitions is described in detail by Knuth in \cite[Section 7.2.1.4]{Knu11} as well as Kelleher \cite{Kell14}. Similarly, Equation \eqref{dioph2} represents the condition that \textbf{m} be a partition of $2de+1$ into squares. This problem (see \cite{OEISseq}) is not as well-studied, although some analyses exist e.g. Bohm et. al. in \cite{Bohm79}. Their simulation experiments suggest that the number of such partitions grows exponentially.

%\begin{figure}
\RestyleAlgo{ruled}
\begin{algorithm}[h]
\caption{\texttt{Extend}: a recursive algorithm for solutions to Equations (\ref{dioph1}) and (\ref{dioph2})}\label{extendalg}
\begin{multicols}{2}
\KwData{$\ell^2$ = remaining $\ell^2$-norm, $\ell^1$ = remaining $\ell^1$-norm, $p$ = previous term}
\KwResult{List of tuples \textbf{Extensions}=\{$\textbf{m}_k$\} }
$\textbf{Extensions} \gets \{\}$\;
%$X \gets x$\;
%$N \gets n$\;
%\DontPrintSemicolon

% \If{$\ell^2$ $\,<\,$ $\ell^1$ \Condition{\Or $n>N$ \;\Or  $\ell^2$ $<0$ \;\Or $\ell^1$ $<0$}}
\If{$\ell^2 < \ell^1$ \textbf{OR} $\ell^2 < 0$ \textbf{OR} $\ell^1<0$}
  {
    \Return \textbf{Extensions}
  }
% \If{$\ell^2$=$\ell^1$ \Condition{\And $N-n>$\,$\ell^2$}}
\If{$\ell^2$=$\ell^1$}
   {
     \Return $\{(1^{\times \ell^2})\}$
   }

\If{$p\neq 0$}
  {
    start = min($p$, $\sqrt{\ell^2})$
  }
\Else{ start $= \sqrt{\ell^2}$ \;
}
  \For{$i\gets$ {\upshape start} \KwTo $1$ \KwBy {-1}}{
     \ForEach{ \textbf{m}$_k$ in \textbf{\upshape Extend}($\ell^2-i^2$, $\ell^1-i$)}{Append \{(i,\textbf{m}$_k$)\} to \textbf{Extensions}}
    }
    %   for i in range(start, 1, -1):
      %       ret = ret + tuple([tuple((i,)+sub) for sub in Extend(square_sum-i*i,linear_sum-i,i, index+1,z)])
      %       if limit!=0:
      %           ret = tuple([x for x in ret if length_bound(x,z)])
      %   return ret
\columnbreak
%\includegraphics[scale=\pagewidth]{}
%\nonl\includestandalone[width=.5\textwidth]{Extend}
\begin{tikzpicture}[%grow=right,
                    level distance=8em,
                    sibling distance=10em,
                    every node/.style = {shape=rectangle, rounded corners,
                    draw,
                    align=center,
                    top color=white
                    %,bottom color=blue!20
                    }]]
\node{Extend($\ell^2=9$, $\ell^1=7$, $p=0$)}
	child { node (I3) {Extend($\ell^2=0$, $\ell^1=4$, $p=3$)}
		child { node (I33) {$(3)$} }
		%child { node {$()$}}
		%child { node {$\frac{1}{3}<d\leq \frac{1}{2}$}
		%	child { node {$m>\lambda-2d'$} }
		%	child { node {$m<\lambda-2d'$} } }
		edge from parent node[left, draw=none, xshift=-3pt] {$i=3$}
        }
	child { node (I2) at ($(I33)+(4,0)$) {Extend($\ell^2=5$, $\ell^1=5$, $p=2$)}
            child {node {$(2,1^{\times 5})$}}
            edge from parent node[right, draw=none, xshift=2pt] {$i=2$}
            };
  \end{tikzpicture}

\nonl Example recursion tree, producing the class $(2,2; 2, 1^{\times 5})$.
\end{multicols}
\end{algorithm}

In both cases, the problem exhibits a recursive structure which we outline in our setting. For fixed $(d,e)$ the self-intersection and Chern number conditions in equations (\ref{dioph1}) and (\ref{dioph2}) constrain the possible $\textbf{m}$ which can appear. For example, by (\ref{dioph2}), the first term $m_1$ in \textbf{m} can be no larger than $\lfloor\sqrt{2de+1}\rfloor$. If $m_1$ takes this maximal value, the subsequent term can be no larger than $\left\lfloor\sqrt{2de+1-\lfloor\sqrt{2de+1}\rfloor^2}\right\rfloor$.

This recursive structure suggests a straightforward method for enumerating candidate obstructive classes. Given a target $\ell^1$ and $\ell^2$ norm for a vector \textbf{m} coming from the Diophantine equations (\ref{dioph1}) and (\ref{dioph2}), we find the largest possible first term in \textbf{m}, i.e. the $\lfloor\sqrt{2de+1}\rfloor$ mentioned above. We then run \texttt{Extend(}$\ell^2-i^2, \ell^1-i)$ for $1\leq i \leq \lfloor\sqrt{2de+1}\rfloor$ to test all possible subsequent entries of \textbf{m}. We also keep track of the previous term to ensure that \textbf{m} is a non-increasing sequence.

Pseudocode for an algorithm \texttt{Extend()} carrying this out appears in Algorithm \ref{extendalg}, along with an example recursion tree. Computer code (written in Python) used to perform this search is available on  \href{https://github.com/crispfish/symplectic-embeddings-cubes/}{GitHub}. While this method will clearly produce all possible classes eventually, it has the downside of a high time complexity.

%Indeed, we may estimate an upper bound for this runtime in terms of $e$. Using $d\leq 2e+2$ and the right-hand sides of (\ref{dioph1}) and (\ref{dioph2}), the upper bound for the $\ell^2$ norm as an input of \texttt{Extend()} is $2e^2+1$. As the loop in line 15 of (\ref{extendalg}) may call the function $\lfloor\sqrt{2e^2+1}\rfloor$ times, this worst-case time complexity is $\mathcal{O}(e^{2e^2+1})$.

%Fortunately, there are improvements possible on the version described above and shown in Algorithm \ref{extendalg}. Lemma \ref{lengthbound} aids in computation by obviating the need to calculate past the bound on $q$ from Lemma \ref{errorestimates}, of which there are finitely many using Lemma \ref{ebound}. Indeed, (\ref{ebound}) bounds $q$ by 11, and the weight expansion of any number of the form $8\frac{n}{11}$ has length at most 19, so there is no need to to continue a search past this point. Even so, given
After generating all possibilities of $\langle d,e;m_1,\dots\rangle$ for a given $d$ and $e$, we use Lemma \ref{oneblock} to eliminate those with more than one block on which the entries differ by more than 1, in accordance with Lemma \ref{oneblock}. As this depends on the weight sequence of $a$, we verify this for each of the possible $a$-values. Even with the bounds $q\leq 11$ and $e\leq 40$, the time and memory demands of our overall method necessitate high-performance computing resources both to generate possible classes and check them for possible obstructions.

Calculating a precise upper bound on the time complexity of \texttt{Extend()} is outside the scope of this paper, but the results of an experiment comparing runtimes for various algorithms appears on \href{https://github.com/crispfish/symplectic-embeddings-cubes/}{GitHub} in the notebook \texttt{timecomparison.ipynb}. The experiment included three algorithms implemented in Python: Hindenburg's algorithm as described in \cite[Section 7.2.1.4]{Knu11}, \texttt{ACCELGEN} from \cite{Kell14}, and \texttt{Extend()} as described above. Given target $\ell^2$ and $\ell^1$ norms for candidate vectors \textbf{m}, \texttt{Extend()} produced results on average at least as fast than both algorithms based on integer partitioning. For this reason, we structured the rest of our search around this method. We note that further optimizations of this algorithm are likely possible, but also out of scope for the current work.

\section{Embedding capacity equals volume for \texorpdfstring{$a\geq 9$}{a greater than 9}}

In this section we use repeated Cremona transforms (along with assumptions or known facts about the ordering of terms) to show that the volume constraint is sufficient to guarantee embeddings if $a \ge 9$ and $1 \le b \le 2.$ Specifically, we prove the following result:

\begin{thm}\label{agreaterthan9}
\textit{If $b \in [1, 2],$ and $a \in [9, +\infty),$ then $c_b(a) = \sqrt{\frac{a}{2b}}.$}
\end{thm}

% \begin{defi}
% A \textit{Cremona move} on an element of $(d;m_1, \dots, m_n)\in H^2(X_t; \mathbb{R})$ is a Cremona transform followed by any permutation of the terms $m_1, \dots, m_n$ of the tail. A \textit{standard Cremona move} is a Cremona transform followed by a permutation of the $m_i$ such that $m_1\geq m_2\geq \ldots \geq m_n$. 
% \end{defi}

We follow the notational convention in \cite{JL} and write \be C_k = (h; t_1, t_2 || t_3, t_4, \ldots)\ee if $t_1\geq t_2$ and both are known to be larger than all of the terms to the right of the symbol $||$. If $C_k$ is an ordered weight vector, we write \be C_k = (h; t_1, t_2, t_3, \ldots),\ee and we let \be \delta_k = h - t_1 - t_2 - t_3 \ee be the defect of $C_k.$

This notation will clarify that we apply Cremona transforms only for vectors where we know the three largest terms $t_1, t_2, t_3$ of the tail. As it is only necessary to know these three largest terms to apply Theorem \ref{reductionvector}, this information (along with verifying that the terms are positive) suffices for guaranteeing an embedding.

\subsection{Proof of Theorem \ref{agreaterthan9}}

Let $a \ge 9,$ and rational; as $c_b(a)$ is continuous it suffices to consider rational $a$. Let \be\w(a) = (1^{\times \floor{a}}, w_1^{\times \ell_1}, \ldots, w_k^{\times \ell_k})\ee be the weight expansion of $a.$ Let $1 \le b \le 2.$

Note that
\begin{equation}
    w_1 = a - \floor{a}.
\end{equation}

Recall from Equation \ref{reductionvector} of Theorem \ref{reductionThm} that the relevant class for the embedding question $E(1,a)\hookrightarrow \lambda P(1,b)$ is of the form
\begin{equation} C_0 = ((b+1)\lambda || b\lambda, \lambda, \w(a)) = ((b+1)\lambda || b\lambda, \lambda, 1^{\times \floor{a}}, w_1^{\times \ell_1}, \ldots, w_t^{\times \ell_t}), \end{equation} where

\begin{equation}
\lambda = \sqrt{\frac{a}{2b}}.
\end{equation}

Now $b\lambda \ge \lambda$ because $b \ge 1.$ Since $a \ge 9,$ we have $b\lambda > 2,$ so it follows that \begin{equation} b\lambda - 1 > 1.\end{equation} Either $\lambda \ge 2$ or not. If $\lambda \ge 2,$ then of course \be C_0 = ((b+1)\lambda; b\lambda, \lambda, 1^{\times \floor{a}}, \ldots),\ee with defect $\delta_0 = -1,$ and so after performing a standard Cremona move, we obtain \be C_1 = ((b+1)\lambda - 1; b\lambda - 1, \lambda - 1, 1^{\times \floor{a} - 1}, \ldots),\ee with defect $\delta_1 = 0.$ Note that if $a \ge 16,$ then $a \ge 8b$ so $\lambda \ge 2.$ In other words, if $a \ge 16,$ then $c_b(a) = \lambda$ and we are done. Thus, we may assume $\lambda < 2$ (hence $9 \le a < 16$).

Under this assumption, we have \be b\lambda - 1 > 1 > \lambda - 1,\ee and \be C_1 = ((b+1)\lambda - 1|| b\lambda - 1, 1^{\times \floor{a} - 1}, \lambda - 1, w_1^{\times \ell_1}, \ldots),\ee where the only remaining ambiguity in the ordering of the first several terms of $C_1$ is whether $\lambda - 1 \ge w_1$ (Case I) or $\lambda - 1 < w_1$ (Case II).

Before handling the cases, we establish a few useful facts that will make our computations easier and more straightforward. First, it will be useful to know when, for positive real numbers $p, q,$ we have $(b + p)\lambda - q \ge 0$ for all $b \in [1, 2].$ 

Now, $(b + p)\lambda - q \ge 0 \iff a \ge 2q^2b/(b + p)^2.$ A straightforward analysis of the function $f(x) = 2q^2x/(x + p)^2$ on the interval $[1, 2]$ shows that it has an absolute maximum at $x = 2$ if $p \ge 2.$

So \be (b + p) \lambda - q \ge 0 \iff a \ge 4q^2/(2 + p)^2.\ee

In light of the above inequalities, if $p$ and $q$ are positive real numbers with $p \ge 2,$ define
\be\label{mpq}
m(p, q) := \frac{4q^2}{(2 + p)^2}.
\ee
In summary, we have:
\begin{lem}\label{mPQ}\textit{
With $m(p,q)$ defined above for positive real numbers $p, q,$ with $p \ge 2,$ we have $(b+p)\lambda -q \ge 0$ if and only if $a \ge m(p,q).$}
\end{lem}

The next lemma will help us determine the correct initial order of various $C_k:$
\begin{lem}\label{biggerThanHalf}\textit{We have $\lambda - 1 > 1/2.$ In particular, if $a$ is such that $w_1 > \lambda - 1,$ and $\ell_1$ is the multiplicity of $w_1$ in the weight expansion of $a,$ then $\ell_1 = 1.$}  \end{lem}
\begin{proof}
    We have $\lambda > 3/2 \iff 2a > 9b.$ Since $a \ge 9$ and $b \le 2,$ we have $2a \ge 18 \ge 9b.$
\end{proof}

For each $k \in \mathbb Z,$ define \be h(k) = b\lambda - 1 + k(\lambda - 2).\ee

\begin{lem}\label{Hkstuff}\textit{For each $k \in \mathbb Z,$ we have:
\begin{enumerate}
\item $h(k+1) = h(k) + \lambda - 2$ and $h(k) \ge 1 \implies h(k + 1) \ge \lambda - 1.$
\item $h(k) + 3 - 2\lambda  = h(k - 2) - 1.$
\item \textit{Suppose $\lambda < 2$ and $k \ge 0$ is the largest integer such that both $h(k) \ge 1$ and $\lfloor a \rfloor - (2k + 1) \ge 0.$  If $\lambda - 1 \ge w_1,$ then \be C_{k+1} = (h(k) + \lambda; h(k), 1^{\times \floor{a} - (2k + 1)}, (\lambda - 1)^{\times 2k + 1}, w_1^{\times \ell_1} \ldots),\ee and if $w_1 > \lambda - 1,$ then $\ell_1 = 1$ and \be C_{k+1} = (h(k) + \lambda; h(k), 1^{\times \floor{a} - (2k + 1)}, w_1, (\lambda - 1)^{\times 2k + 1}, \ldots).\ee}
\end{enumerate}}
\end{lem}

\begin{proof} Part (1) is straightforward. For (2), we have
\begin{eqnarray*}
h(k) + 3 - 2\lambda &=& h(k) + (4 - 2\lambda) - 1\\
&=& (h(k) + (2 - \lambda)) + (2 - \lambda) - 1\\
&=& (h(k- 1) + 2 - \lambda) - 1\\
&=& h(k-2) - 1.
\end{eqnarray*} 
To see (3), first note that \be C_1 = (h(0) + \lambda; h(0), 1^{\times \lfloor a \rfloor -1}||(\lambda - 1)^{\times 1}, w_1^{\times \ell_1},\ldots),\ee with defect $\delta_1 = \lambda - 2.$ Suppose \be C_{j+1} = (h(j) + \lambda; h(j), 1^{\lfloor a \rfloor - (2j + 1)}|| (\lambda - 1)^{\times (2j + 1)}, w_1^{\ell_1}, \ldots)\ee has been shown for some $0 \le j < k.$ Then the defect of $C_{j+1}$ is $\delta_{j+1} = \lambda - 2,$ and $h(j+1) = h(j) + \lambda - 2$ by part (1), and since $j + 1 \le k,$ we have $h(j+1) \ge h(k) \ge 1.$ 
Therefore, \be
C_{j+2} = C_{(j+1) + 1} = (h(j+1) + \lambda; h(j+1), 1^{\lfloor a \rfloor - (2(j+1) + 1)}|| (\lambda - 1)^{\times (2j + 1)}, w_1^{\ell_1}, \ldots)
\ee
By induction, we have \be C_{k+1} = (h(k) + \lambda; h(k), 1^{\lfloor a \rfloor - (2k+1)}|| (\lambda - 1)^{\times (2k + 1)}, w_1^{\ell_1}, \ldots). \ee

If $\lambda - 1 \ge w_1,$ then the assertion on the order of $C_{k+1}$ is immediate. If $w_1 > \lambda - 1,$ then by Lemma \ref{biggerThanHalf}, $w_1 > 1/2$ (hence $\ell_1 = 1$), and thus $w_2 = 1 - w_1  < 1/2 < \lambda - 1.$ In other words, $\lambda - 1$ bigger than all the remaining terms in the weight expansion of $a.$ So \be C_{k+1} = (h(k) + \lambda; h(k), 1^{\lfloor a \rfloor - (2k+1)}, w_1, (\lambda - 1)^{2k + 1}, \ldots).\ee Finally, note that $h(k) \ge 1 \ge \max \{w_1, \lambda - 1\},$ so if $\lfloor a \rfloor = 2k + 1,$ then the claims about the ordering of $C_{k+1}$ still hold in both cases.
\end{proof}

% \end{proof}

%We begin this case with the following lemma.

We are now ready to handle the cases found above depending on whether $\lambda - 1 \ge w_1$ or $\lambda - 1 < w_1.$ Each of these cases is split up into several subcases where we restrict $a$ to be in a certain interval $I \subseteq [9, 16)$ in order to have greater control over the weight expansion of $a.$ For instance, if we know $a \in [k, k+1)$ for some $k \ge 9,$ then we know $\lfloor a \rfloor = k,$ which is precisely the multiplicity of 1 in the relevant class. The intervals that we choose are influenced by this observation and Lemma \ref{Hkstuff}. In each subcase, we use Lemma \ref{Hkstuff} (3) to reduce the multiplicity of 1 as much as immediately possible. The other identities and inequalities in Lemma \ref{Hkstuff} are used extensively and sometimes without direct reference to the lemma.

\subsection{Case I: \texorpdfstring{$\lambda - 1 \ge w_1$}{lambda minus 1 greater than d1}} In this case, we have \be\label{caseIineq} b\lambda - 1 > 1 > \lambda - 1 \ge w_1\ee and thus \be C_1 = ((b+1)\lambda - 1; b\lambda - 1, 1^{\times \floor{a} - 1}, \lambda - 1, w_1^{\times \ell_1}, \ldots ),\ee with defect $\delta_1 = \lambda - 2.$

\subsubsection{\texorpdfstring{$a \in [9, 10)$}{a between 9 and 10}} We have $h(2) \ge 1$ by Lemma \ref{mPQ} because $a \ge 9 = m(2, 6).$ Applying Lemma \ref{Hkstuff} (3) to $k = 2,$ we have
\be C_3 = (h(2) + \lambda; h(2), 1^{\times 4}, (\lambda - 1)^{\times 5}, w_1^{\times \ell_1}, \ldots),\ee with $\delta_3 = \lambda - 2.$

So \be C_4 = (h(3) + \lambda|| h(3), 1^{\times 2}, (\lambda - 1)^{\times 7}, w_1^{\times \ell_1} \ldots).\ee

If $h(3) \ge 1,$ then $\delta_4 = \lambda - 2,$ and \be C_5 = (h(4) + \lambda;  h(4), (\lambda - 1)^{\times 9}, w_1^{\times \ell_1}, \ldots)\ee because $h(3) \ge 1 \implies h(4) \ge \lambda - 1.$ Therefore, $\delta_5  = 2 - \lambda > 0.$

Now suppose $h(3) < 1.$ Now $h(2) \ge 1 \implies h(3) \ge \lambda - 1,$ and \be 1 > h(3) \implies \lambda - 1 > h(4).\ee Putting everything together, we get $1 > h(3) \ge \lambda - 1 > h(4),$ and \be C_4 = (h(3) + \lambda; 1^{\times 2}, h(3), (\lambda - 1)^{\times 7}, w_1^{\times \ell_1}, \ldots),\ee with $\delta_4 = \lambda - 2.$ Thus \be C_5 = (h(4)+\lambda || (\lambda - 1)^{\times 9}, h(4), w_1^{\times \ell_1}, \ldots).\ee Since $\lambda - 1$ exceeds $h(4)$ and $w_1, \ldots, w_t,$ it follows that \be\delta_5 = h(4) + 3 - 2\lambda = h(2) - 1 \ge 0,\ee where the last equality follows from Lemma \ref{Hkstuff}(2).

\subsubsection{\texorpdfstring{$a \in [10, 10.24)$}{a between 10 and 10.24}} In this situation, $h(2) > 1,$ and by Lemma \ref{Hkstuff} (3) we get \be C_3 = (h(2) + \lambda; h(2), 1^{\times 5}, (\lambda - 1)^{\times 5}, w_1^{\times \ell_1}, \ldots)\ee with $\delta_3 = \lambda - 2,$ so \be C_4 = (h(3) + \lambda|| h(3), 1^{\times 3}, (\lambda - 1)^{\times 7}, w_1^{\times \ell_1}, \ldots).\ee

If $h(3) \ge 1,$ then $\delta_3 = \lambda - 2$ and $h(4) \ge \lambda - 1,$ so \be C_5 = (h(4) + \lambda || h(4), 1^{\times 1}, (\lambda - 1)^{\times 9}, w_1^{\times \ell_1}, \ldots),\ee and thus $\delta_5 = 0.$

Otherwise, $1 > h(3),$ so \be C_4 = (h(3) + \lambda; 1^{\times 3}, h(3), (\lambda - 1)^{\times 7}, w_1^{\times \ell_1}, \ldots),\ee and $\delta_4 = h(3) + \lambda - 3 = h(4) - 1.$ Since $\lambda - 1 > h(4),$ we have $h(4) - 1 <0.$ Therefore, \be C_5 = (h(3) + h(4) + \lambda - 1|| h(3), (\lambda - 1)^{\times 7}, h(4)^{\times 3}, w_1^{\times \ell_1}, \ldots).\ee Since $h(3) > \lambda - 1$ (because $h(2) > 1$), we have $\delta_5 = h(4) + 1 - \lambda  = h(3) -  1 < 0.$ Therefore,

\begin{eqnarray*}
C_6 &=& (h(3) + 2h(4)|| h(3) + h(4) + 1 - \lambda, (\lambda - 1)^{\times 5}, h(4)^{\times 5}, w_1^{\times \ell_1}, \ldots)\\
&=& (3h(3) + 2\lambda - 4|| 2h(3) - 1, (\lambda - 1)^{\times 5}, h(4)^{\times 5}, w_1^{\times \ell_1}, \ldots).
\end{eqnarray*}

Now \be 2h(3) - 1 \ge \lambda - 1 \iff (b + 5/2)\lambda - 7 \ge 0 \iff a \ge m(5/2, 7) \approx 9.68.\ee

Therefore, $\delta_6 = h(3)  - 1 <0,$ so

\be C_7 = (4h(3) + 2\lambda - 5|| 3h(3) - 2, (\lambda - 1)^{\times 3}, h(4)^{\times 7}, w_1^{\times \ell_1}, \ldots).\ee

Similar to above reasoning, we have \be 3h(3) - 2 \ge \lambda - 1 \iff (b + 8/3)\lambda - 22/3 \ge 0 \iff a \ge m(8/3, 22/3) \approx 9.88.\ee

Again, we have $\delta_7 = h(3) - 1,$ so \be C_8 = (5h(3) + 2\lambda - 6|| 4h(3) - 3, (\lambda - 1)^{\times 1}, h(4)^{\times 9}, w_1^{\times \ell_1}, \ldots ).\ee Now \be 4h(3) - 3 \ge \lambda - 1 \iff (b + 11/4)\lambda - 15/2 \ge 0 \iff a \ge m(11/4, 15/2) \approx 9.97.\ee Moreover, $w_1 = a - \floor{a} \le .24,$ and we have \be h(4) \ge .24 \iff (b+4)\lambda - 9.24 \ge 0 \iff a \ge m(4, 9.24) \approx 9.49.\ee

In summary, we have $$4h(3) - 3 \ge \lambda - 1 > h(4) \ge w_1,$$ so
\begin{eqnarray*}
\delta_8 &=& 5h(3) + 2\lambda - 6 - 4h(3) + 3 - \lambda + 1 - h(4)\\
&=& h(3) + \lambda - 2 - h(4)\\
&=& h(4) - h(4) = 0.
\end{eqnarray*}

\subsubsection{\texorpdfstring{$a \in [10.24, 12)$}{a between 10.24 and 12}} First, we note that \be h(3) = (b+3)\lambda - 7 > 1\ee for all $1 \le b \le 2$ because $a \ge m(3, 8) = 10.24.$

Assume $10.24 \le a < 11.$ By Lemma \ref{Hkstuff} (3), we have \be C_4 = (h(3) + \lambda; h(3), 1^{\times 3}, (\lambda - 1)^{\times 7}, w_1^{\times \ell_1}, \ldots)\ee with $\delta_4 = \lambda - 2.$ Therefore, \be C_5  = (h(4) + \lambda || h(4), 1^{\times 1}, (\lambda - 1)^{\times 9}, w_1^{\times \ell_1}, \ldots).\ee Now $h(4) > \lambda - 1$ by Lemma \ref{Hkstuff}(1). In particular, \be\lambda - 1 = \min\{1, \lambda - 1, h(4)\}.\ee Therefore, $\delta_5 = 0.$

Now assume $11 \le a < 12.$ Lemma \ref{Hkstuff} (3) shows \be C_4 = (h(3) + \lambda; h(3), 1^{\times 4}, (\lambda - 1)^{\times 7}, w_1^{\times \ell_1}, \ldots),\ee with defect $\delta_4 = \lambda - 2.$ So \be C_5 = (h(4) + \lambda|| h(4), 1^{\times 2}, (\lambda - 1)^{\times 9}, w_1^{\times \ell_1}, \ldots).\ee

Assume $h(4) \ge 1.$ Recall that $h(4) \ge 1 \implies h(5) \ge \lambda - 1.$ Therefore, \be C_5 = (h(4) + \lambda; h(4), 1^{\times 2}, (\lambda - 1)^{\times 9}, w_1^{\times \ell_1}, \ldots ),\ee with $\delta_5 = \lambda - 2.$ Applying another Cremona move, we see that \be C_6 = (h(5) + \lambda; h(5), (\lambda - 1)^{\times 11}, w_1^{\times \ell_1}, \ldots),\ee with defect $\delta_6 = 2 - \lambda > 0.$

Now suppose $h(4) < 1.$ Similar to previous reasoning, we have
\begin{equation}\label{useful} 1 > h(4) > \lambda - 1 > h(5).\end{equation}

So \be C_5 = (h(4) + \lambda; 1^{\times 2}, h(4), (\lambda - 1)^{\times 9}, w_1^{\times \ell_1}, \ldots),\ee with $\delta_5 = \lambda - 2.$ By \ref{useful}, we get \be C_6 = (h(5) + \lambda|| (\lambda - 1)^{\times 11}, h(5), w_1^{\times \ell_1}, \ldots),\ee with $\delta_6 = h(5) + 3 - 2\lambda = h(3) - 1 > 0,$ where the latter equality follows from Lemma \ref{Hkstuff}(2).

\subsubsection{\texorpdfstring{$a \in [12, 16)$}{a between 12 and 16}} We start by making the following observations:

\begin{eqnarray*}
h(5) \ge 1 &\iff& a \ge m(5, 12) = 11.755.\\
h(6) \ge 1 &\iff& a \ge m(6,14) = 12.25.\\
h(7) \ge 1 &\iff& a \ge m(7, 16) \approx 12.642.\\
%h(8) \ge 1 &\iff& a \le m(8, 18) = 12.96.
\end{eqnarray*}

If $12 \le a < 13,$ then $h(5) \ge 1,$ so Lemma \ref{Hkstuff} (3) implies \be C_6 = (h(5) + \lambda; h(5), 1^{\times 1}, (\lambda - 1)^{\times 11}, w_1^{\times \ell_1}, \ldots)\ee with $\delta_6 = 0.$

If $13 \le a < 14,$ then $h(6) \ge 1,$ so \be C_7= (h(6) + \lambda; h(6), (\lambda - 1)^{\times 13}, w_1^{\times \ell_1}, \ldots),\ee with $\delta_7 = 2 - \lambda>0.$

For $14 \le a < 15,$ we still use $h(6) \ge 1$ and Lemma \ref{Hkstuff} (3) to get \be C_7= (h(6) + \lambda; h(6), 1^{\times 1}, (\lambda - 1)^{\times 13}, w_1^{\times \ell_1}, \ldots),\ee with $\delta_7 = 0.$

Finally, if $15 \le a < 16,$ then $h(7) \ge 1$ and Lemma \ref{Hkstuff}(3) gives \be C_8 = (h(7) + \lambda; h(7), 1^{\times 1}, (\lambda -1)^{\times 14}, w_1^{\times \ell_1}, \ldots),\ee with $\delta_8 = 0.$

\subsection{Case II: \texorpdfstring{$w_1 > \lambda - 1$}{d1 greater than lambda minus 1}}

\subsubsection{\texorpdfstring{$a \in [9, 10)$}{a between 9 and 10}}

Since $h(2) \ge 1,$ Lemma \ref{Hkstuff}(3) gives us \be C_3 = (h(2) + \lambda; h(2), 1^{\times 4}, w_1, (\lambda - 1)^{\times 5}, \ldots)\ee with defect $\delta_3 = \lambda - 2 < 0.$ So \be C_4 = (h(3) + \lambda || h(3), 1^{\times 2}, w_1, (\lambda - 1)^{\times 7}, \ldots).\ee

From this point, there are several cases to consider, each of which yields different calculations. The cases stem from how $h(3), 1,$ and $w_1$ relate to one another. Each case below is viewed as a subcase of the assumption that $a \in [9, 10).$ Such will be the case in future subsections where similar computations are demonstrated. Computations of subsequent $\delta_k$ will make liberal use of the identities in Lemma \ref{Hkstuff}.

\textbf{Subcase 1: $h(3) \ge 1.$} This case is the most straightforward. In this case, \be h(3) \ge 1 \ge w_1 > \lambda - 1\ee so $\delta_4 = \lambda - 2 < 0.$ Therefore, \be C_5 = (h(4) + \lambda || h(4), w_1, (\lambda - 1)^{\times 9}, \ldots).\ee
Since $h(3) \ge 1,$ we have $h(4) \ge \lambda - 1$ by Lemma \ref{Hkstuff}(1). Thus $\delta_5 = 1 - w_1 \ge 0.$ %Otherwise, if $\lambda - 1 > h(4),$ then $w_1 > \lambda - 1 > h(4),$ so
%\begin{eqnarray*}
%\delta_5 &=& h(4) + \lambda - w_1 -2(\lambda - 1)\\
%&=& h(4) - w_1 + 2 - \lambda \\
%&\ge& h(4) - 1 + 2 - \lambda \\
%&=& h(3) - 1 \ge 0.
%\end{eqnarray*}

\textbf{Subcase 2: $1 > h(3) > w_1.$} This case is very similar to the previous since, as before, we have $\delta_4 = \lambda - 2,$ and  \be C_5 = (h(4) + \lambda || h(4), w_1, (\lambda - 1)^{\times 9}, \ldots).\ee The analysis directly above shows that $\delta_5 \ge 0$ in this case.

\textbf{Subcase 3: $1 > w_1 > h(3).$}
To begin, we need a key inequality.

\begin{lem}\label{910inequalities} If $w_1 > h(3)$ and $9 \le a < 10,$ we have \begin{equation}\label{910inequality}
    4h(3) - 3w_1 > \lambda - 1 > h(4) + 1 - w_1.
\end{equation}
\end{lem}

\begin{proof}
    On $[9, 10),$ $w_1 = a - 9,$ so the first inequality in \ref{910inequality} holds if and only if \begin{equation}\label{ineq1}
    (4b+11)\lambda - 3a > 0. \end{equation} Fix $b,$ and view the left-hand side of \ref{ineq1} as a function of $a.$ The derivative of the left-hand side is $(4b+11)\lambda'-3,$ which is negative if and only if $\lambda' = \frac{1}{4b}\sqrt{\frac{2b}{a}} < 3/(4b+11),$ or equivalently, \be\sqrt{\frac{2b}{a}} < \frac{12b}{4b + 11}.\ee Now $12b/(4b+11) \ge 12/15$ for $b \in [1, 2],$ and $\sqrt{2b/a} \le \sqrt{2b}/3$ on $[9, 10).$ The above inequality will hold, then, if $\sqrt{2b} < 36/15,$ which is certainly true. Therefore, $(4b + 11)\lambda - 3a$ is decreasing as a function of $a$ on $[9, 10).$ Now \be\lim_{a \to 10^{-}} (4b + 11)\lambda - 3a = (4b + 11)\sqrt{5/b} - 30 > 0,\ee where the latter inequality follows from familiar considerations using differential calculus. In summary, the function in question is positive on $[9, 10),$ and the first inequality in \ref{910inequality} is hence established.

    For the second inequality in \ref{910inequality}, note that $\lambda - 1 > h(4) + 1 - w_1$ if and only if $w_1 > h(4) + 2 - \lambda = h(3),$ which is given. \end{proof}

As an immediate consequence of Lemma \ref{910inequalities} and our key assumption that $w_1 > h(3),$ we have \begin{equation}\label{bigIneq}
    h(3) > 2h(3) - w_1 > 3h(3) - 2w_1 > 4h(3) - 3w_1 > \lambda - 1 > h(4) + 1 - w_1.
\end{equation}

With these inequalities in mind, we are ready to begin reduction calculations in this case. Throughout, recall $w_2 = 1 - w_1$ from the weight expansion of $a.$ We have \be C_4 = (h(3) + \lambda; 1^{\times 2}, w_1, h(3), (\lambda - 1)^{\times 7}|| \ldots),\ee so that $\delta_4 = h(3) + \lambda  - 2 - w_1 = h(4) - w_1.$

If $\delta_4 \ge 0,$ we are done so assume $\delta_4 < 0.$ Then \be C_5 = (h(3) + h(4) + \lambda - w_1 || (h(4) + w_2)^{\times 2}, h(4), h(3), (\lambda - 1)^{\times 7}, \ldots).\ee Now $1 > h(3) \implies \lambda - 1 > h(4),$ so \be C_5 = (h(3) + h(4) + \lambda - w_1; h(3), (\lambda - 1)^{\times 7}, (h(4) + w_2)^{\times 2} || \ldots ),\ee and $\delta_5 = h(3) - w_1 < 0.$ Now $\lambda - 1 + h(3) - w_1 = h(4) + w_2$ which, combined with \ref{bigIneq}, implies \be C_6 = (2h(3) + h(4) + \lambda - 2w_1; 2h(3) - w_1, (\lambda - 1)^{\times 5}, (h(4) + w_2)^{\times 4} || \ldots),\ee so that again $\delta_6 = h(3) - w_1 < 0.$ Continuing, we see \be C_7 = (3h(3) + h(4) + \lambda - 3w_1; 3h(3) - 2w_1, (\lambda - 1)^{\times 3}, (h(4) + w_2)^{\times 6}||\ldots),\ee so that $\delta_7 = h(3) - w_1 < 0.$ Finally, we see that \be C_8 = (4h(3) + h(4) + \lambda - 4w_1; 4h(3) - 3w_1, (\lambda - 1)^{\times 1}, (h(4) + w_2)^{\times 8} || \ldots),\ee and $\delta_8 = 0.$ This concludes the case where $9 \le a < 10.$

\subsubsection{\texorpdfstring{$a \in [10, 11)$}{a between 10 and 11}} In this case, we claim that $a > 10.5$ given $w_1 > \lambda - 1.$ Indeed, $w_1 = a - 10$ and since $\lambda - 1 > 1/2$ by Lemma \ref{biggerThanHalf}, we have $a - 10 > 1/2.$ In particular, $h(3) \ge 1$ because $m(3,8) = 10.24 < 10.5.$ By Lemma \ref{Hkstuff}(1), we have $h(4) \ge \lambda - 1$ and \be C_4 = (h(3) + \lambda; h(3), 1^{\times 3}, w_1, (\lambda - 1)^{\times 7}|| \ldots),\ee with defect $\delta_4 =  \lambda - 2$ so that \be C_5 = (h(4) + \lambda||h(4), 1^{\times 1}, w_1, (\lambda -1 )^{\times 9}, \ldots).\ee

After re-arranging the first three terms past the head, we get \be C_5 = (h(4) + \lambda|| 1^{\times 1}, w_1, h(4), (\lambda -1 )^{\times 9}, \ldots).\ee Since $h(4) \ge \lambda - 1$ and $w_1 > \lambda - 1,$ we have $\min\{1, w_1, h(4)\} \ge \lambda - 1,$ so $\delta_5 = \lambda - 1 - w_1 < 0,$ and we get \be C_6 = (h(5) + \lambda + 1 - w_1 || h(5) + 1 - w_1, \lambda - w_1, (\lambda - 1)^{\times 10}, \ldots).\ee

If $h(5) + 1 - w_1 \ge \lambda - 1,$ then $\delta_6 = w_1 - \lambda + 1 >0.$ Otherwise, $h(5) + 1 - w_1 < \lambda - 1,$ and we have \be\delta_6 = h(5) + 3 - 2\lambda  = h(4) + 1 - \lambda \ge 0.\ee

\subsubsection{\texorpdfstring{$a \in [11, 12)$}{a between 11 and 12}} Before we get started here, we claim $h(5) > w_1.$ This is equivalent to showing $a > m(5, 11+w_1).$ Since $w_1 = a - 11,$ we have \be m(5, 11+w_1) = \frac{4(11 + w_1)^2}{7^2} = \frac{4}{49}a^2.\ee So $a > m(5, 11+w_1)$ if and only if $a < 49/4 = 12.25.$ So $h(5) > w_1$ hence $h(4) > h(5) > w_1$ as well. 

Since $h(3) \ge 1,$ we have by Lemma \ref{Hkstuff}(3) \be C_4 = (h(3) + \lambda; h(3), 1^{\times 4}, w_1, (\lambda - 1)^{\times 7}|| \ldots ),\ee with $\delta_4 = \lambda - 2.$ So, \be C_5= (h(4) + \lambda|| h(4), 1^{\times 2}, w_1, (\lambda - 1)^{\times 9} \ldots).\ee

Since $\min\{h(4), 1\} \ge w_1,$ we have $\delta_5 = \lambda - 2$ and

\be C_6 = (h(5) + \lambda; h(5), w_1, (\lambda - 1)^{\times 11}|| \ldots),\ee with $\delta_6 = 1 - w_1 \ge 0.$

\subsubsection{\texorpdfstring{$a \in [12, 13)$}{a between 12 and 13}} Since $h(5) \ge 1$ on this interval, we have by Lemma \ref{Hkstuff}(3) \be C_6 = (h(5) + \lambda; h(5), 1^{\times 1}, w_1, (\lambda - 1)^{\times 11}|| \ldots)\ee so $\delta_6 = \lambda - 1 - w_1 < 0.$ Then \be C_7 = (h(6) + \lambda + 1 - w_1; h(6)+ 1- w_1, \lambda - w_1, (\lambda - 1)^{\times 12}|| \ldots),\ee with $\delta_7 = w_1 - \lambda + 1 > 0.$

\subsubsection{\texorpdfstring{$a \in [13, 14)$}{a between 13 and 14}} Since $h(6) \ge 1 \ge w_1,$ we get by Lemma \ref{Hkstuff}(3) \be C_7 = (h(6) + \lambda; h(6), w_1, (\lambda - 1)^{\times 13}||\ldots),\ee so that $\delta_7 = 1 - w_1 \ge 0.$

\subsubsection{\texorpdfstring{$a \in [14, 15)$}{a between 14 and 15}} Again, we have $h(6) \ge 1$ so \be C_7 = (h(6) + \lambda; h(6), 1^{\times 1}, w_1, (\lambda - 1)^{\times 13}|| \ldots).\ee Then $\delta_7 = \lambda - 1 - w_1 < 0$ so that \be C_8 = (h(7) + \lambda + 1 - w_1; h(7) + 1 - w_1, \lambda - w_1, (\lambda - 1)^{\times 14}||\ldots).\ee

In particular, $\delta_8 = w_1 - \lambda + 1 > 0.$

\subsubsection{\texorpdfstring{$a \in [15, 16)$}{a between 15 and 16}} In this case, we have $h(7) \ge 1$ so that \be C_8 = (h(7) + \lambda|| h(7), 1^{\times 1}, w_1, (\lambda - 1)^{14}, \ldots ).\ee

From here, we see $\delta_8 = \lambda - 1 - w_1 < 0,$ and thus \be C_9 = (h(8) + \lambda + 1 - w_1|| h(8) + 1 - w_1, \lambda - w_1, (\lambda - 1)^{\times 15},...),\ee  and thus $\delta_ 9 = w_1 - \lambda + 1 > 0.$

\printbibliography

\end{document}